\documentclass{amsart}

\usepackage{amsmath} 
\usepackage{amssymb}

\newtheorem{theorem}{Theorem}[section] 
\newtheorem{claim}{Claim}[theorem]
 
\newtheorem{proposition}[theorem]{Proposition} 
\newtheorem{observation}[theorem]{Observation} 
\newtheorem{corollary}[theorem]{Corollary} 

\theoremstyle{definition}
\newtheorem{definition}[theorem]{Definition}

\newtheorem{problem}[theorem]{Problem}

\theoremstyle{remark}
\newtheorem{remark}[theorem]{Remark}

\newtheorem{conclusion}[theorem]{Conclusion}
\newtheorem{context}[theorem]{Context}

\numberwithin{equation}{section}
\newcommand{\forces}{\Vdash}

\newcommand{\bV}{{\bf V}} 
\newcommand{\lesdot}{\mathrel{\mathord{<}\!\!\raise 
0.8 pt\hbox{$\scriptstyle\circ$}}}

\newcommand{\conc}{{}^\frown\!}
\newcommand{\lh}{{\rm lh}\/}
\newcommand{\rest}{{\restriction}}
\newcommand{\Dom}{{\rm Dom}}

\newcommand{\gb}{{\mathfrak b}}

\newcommand{\gc}{{\mathfrak c}}

\newcommand{\cD}{{\mathcal D}}
\newcommand{\gd}{{\mathfrak d}} 

\newcommand{\cH}{{\mathcal H}}

\newcommand{\cF}{{\mathcal F}}

\newcommand{\cI}{{\mathcal I}}

\newcommand{\bbP}{{\mathbb P}}
\newcommand{\cP}{{\mathcal P}}

\newcommand{\bbQ}{{\mathbb Q}}

\newcommand{\cT}{{\mathcal T}}
\newcommand{\cU}{{\mathcal U}}
\newcommand{\cV}{{\mathcal V}}

\newcommand{\cl}{{\rm cl}\/}

\newcommand{\otp}{{\rm otp}\/}

\newcommand{\st}{{\bf st}} 

\newcommand{\vare}{\varepsilon}

\newcommand{\tegame}{{\Game^{\rm rcB}_{\cU,\bar{\mu}}}}
\newcommand{\tegamel}{{\Game^{\rm rcB}_{\cU,\bar{\mu}'}}}
\newcommand{\Agame}{{\Game^{\rm rcA}_{\bar{\mu}}}}
\newcommand{\Cgame}{{\Game^{\rm rcC}_{\cU,\bar{\mu}}}}

\newcommand{\agame}{{\Game^{{\rm rc}{\bf a}}_{\bar{\mu}}}}
\newcommand{\bgame}{{\Game^{{\rm rc}{\bf b}}_{\cU,\bar{\mu}}}} 
\newcommand{\cgame}{{\Game^{{\rm rc}{\bf c}}_{\cU,\bar{\mu}}}}
\newcommand{\tagame}{{\Game^{{\rm rc}{\bf 2a}}_{\bar{\mu}}}}
\newcommand{\tbgame}{{\Game^{{\rm rc}{\bf 2b}}_{\cU,\bar{\mu}}}}

\newcommand{\clfo}{{\bbQ^1_S}}
\newcommand{\tefo}{{\bbQ^2_{\cU}}}

\newcommand{\pos}{{\rm pos}}

\newcommand{\mrot}{{\rm root}}

\newcommand{\rk}{{\rm rk}}
\newcommand{\vtl}{\vartriangleleft}

\newcommand{\Gsg}{{\Game^{\rm Sacks}_{\bar{\mu}}}}
\newcommand{\bairel}{{}^{\lambda}\lambda}

\newcount\skewfactor
\def\mathunderaccent#1#2 {\let\theaccent#1\skewfactor#2
\mathpalette\putaccentunder}
\def\putaccentunder#1#2{\oalign{$#1#2$\crcr\hidewidth
\vbox to.2ex{\hbox{$#1\skew\skewfactor\theaccent{}$}\vss}\hidewidth}}
\def\name{\mathunderaccent\tilde-3 }

\begin{document}

\title{Reasonably complete forcing notions}

\author{Andrzej Ros{\l}anowski}
\address{Department of Mathematics\\
 University of Nebraska at Omaha\\
 Omaha, NE 68182-0243, USA}
\email{roslanow@member.ams.org}
\urladdr{http://www.unomaha.edu/logic}

\author{Saharon Shelah}
\address{Einstein Institute of Mathematics\\
Edmond J. Safra Campus, Givat Ram\\
The Hebrew University of Jerusalem\\
Jerusalem, 91904, Israel\\
 and  Department of Mathematics\\
 Rutgers University\\
 New Brunswick, NJ 08854, USA}
\email{shelah@math.huji.ac.il}
\urladdr{http://www.math.rutgers.edu/$\sim$shelah}
\thanks{Both authors acknowledge support from the United States-Israel
Binational Science Foundation (Grant no. 2002323). This is publication
860 of the second author.}

\date{June 2005}

\begin{abstract}
We introduce more properties of forcing notions which imply that
their $\lambda$--support iterations are $\lambda$--proper, where $\lambda$
is an inaccessible cardinal. This paper is a direct continuation of
Ros{\l}anowski and Shelah \cite[\S A.2]{RoSh:777}. As an application of our
iteration result we show that it is consistent that dominating numbers
associated with two normal filters on $\lambda$ are distinct. 
\end{abstract}

\maketitle

\section{Introduction}
There are serious ZFC obstacles to easy generalizations of properness to the
case of iterations with uncountable supports (see, e.g., Shelah
\cite[Appendix 3.6(2)]{Sh:f}). This paper belongs to the series of works
aiming at localizing ``good properness conditions'' for such iterations and
including Shelah \cite{Sh:587}, \cite{Sh:667}, Ros{\l}anowski and Shelah
\cite{RoSh:655}, \cite{RoSh:777} and Eisworth \cite{Ei03}. This paper is a
direct continuation of Ros{\l}anowski and Shelah \cite[\S A.2]{RoSh:777},
though no familiarity with the previous paper is assumed and the current
work is fully self-contained. 

In Section 2 we introduced 3 bounding--type properties (A, B, C) and we
essentially show that the first two are almost preserved in
$\lambda$--support iterations (Theorems \ref{verB}, \ref{verA}). ``Almost''
as the limit of the iteration occur to have somewhat weaker property, but
equally applicable. In the following section we show that the {\em
  reasonable A--bounding\/} property is equivalent to the one introduced in
\cite[\S A.2]{RoSh:777} thus showing that \ref{verA} improves \cite[Theorem
  A.2.4]{RoSh:777}. Finally, in the fourth section of the paper, we give an
example of an interesting {\em reasonable $B$--bounding\/} forcing notion
and we use it  to show that it is consistent that dominating numbers
associated with two normal filters on $\lambda$ are distinct (Conclusion
\ref{conc}).   

Like in \cite{RoSh:777}, we assume here that our cardinal $\lambda$ is
inaccessible. We do not know at the moment if any parallel work can be done
for a successor cardinal. 
\medskip

\noindent {\bf Notation:}\quad Our notation is rather standard and
compatible with that of classical textbooks (like Jech \cite{J}). In forcing
we keep the older convention that {\em a stronger condition is the larger
  one}.  

\begin{enumerate}
\item Ordinal numbers will be denoted be the lower case initial letters of
the Greek alphabet ($\alpha,\beta,\gamma,\delta\ldots$) and also by $i,j$
(with possible sub- and superscripts). 
 
Cardinal numbers will be called $\kappa,\lambda,\mu$; {\em $\lambda$ will be  
always assumed to be inaccessible\/} (we may forget to mention it).

By $\chi$ we will denote a {\em sufficiently large\/} regular cardinal; 
$\cH(\chi)$ is the family of all sets hereditarily of size less than
$\chi$. Moreover, we fix a well ordering $<^*_\chi$ of $\cH(\chi)$. 

\item For two sequences $\eta,\nu$ we write $\nu\vartriangleleft\eta$
whenever $\nu$ is a proper initial segment of $\eta$, and $\nu
\trianglelefteq\eta$ when either $\nu\vartriangleleft\eta$ or $\nu=\eta$. 
The length of a sequence $\eta$ is denoted by $\lh(\eta)$.

\item We will consider several games of two players. One player will be
called {\em Generic\/} or {\em Complete\/} or just {\em COM\/}, and we
will refer to this player as ``she''. Her opponent will be called {\em
Antigeneric\/} or {\em Incomplete} or just {\em INC\/} and will be
referred to as ``he''. 

\item For a forcing notion $\bbP$, $\Gamma_\bbP$ stands for the canonical
$\bbP$--name for the generic filter in $\bbP$. With this one exception, all
$\bbP$--names for objects in the extension via $\bbP$ will be denoted with a
tilde below (e.g., $\name{\tau}$, $\name{X}$). The weakest element of $\bbP$
will be denoted by $\emptyset_\bbP$ (and we will always assume that there is
one, and that there is no other condition equivalent to it). We will also
assume that all forcing notions under considerations are atomless. 

By ``$\lambda$--support iterations'' we mean iterations in which domains of
conditions are of size $\leq\lambda$. However, we will pretend that
conditions in a $\lambda$--support iteration $\bar{\bbQ}=\langle\bbP_\zeta,
\name{\bbQ}_\zeta:\zeta<\zeta^*\rangle$ are total functions on $\zeta^*$
and for $p\in\lim(\bar{\bbQ})$ and $\alpha\in\zeta^*\setminus\Dom(p)$ we
will let $p(\alpha)=\name{\emptyset}_{\name{\bbQ}_\alpha}$. 

\item For a filter $D$ on $\lambda$, the family of all $D$--positive subsets
of $\lambda$ is called $D^+$. (So $A\in D^+$ if and only if $A\subseteq
\lambda$ and $A\cap B\neq\emptyset$ for all $B\in D$.)

The club filter of $\lambda$ is denoted by $\cD_\lambda$.

\end{enumerate}
\bigskip

 \begin{context}
\label{incon}
\begin{enumerate}
\item[(a)] $\lambda$ is a strongly inaccessible cardinal,
\item[(b)] $\bar{\mu}=\langle\mu_\alpha:\alpha<\lambda\rangle$, each
$\mu_\alpha$ is a regular cardinal satisfying (for $\alpha<\lambda$)
\[\aleph_0\leq\mu_\alpha\leq\lambda\qquad\mbox{ and }\qquad \big(\forall f\in
{}^\alpha \mu_\alpha\big)\big(\big|\prod_{\xi<\alpha} f(\xi)\big|<
\mu_\alpha\big),\]
\item[(c)] $\cU$ is a normal filter on $\lambda$. 
\end{enumerate}
\end{context}

\section{Preliminaries on $\lambda$--support iterations}

\begin{definition}
\label{strcom}
Let $\bbP$ be a forcing notion.
\begin{enumerate}
\item For a condition $r\in\bbP$ let $\Game_0^\lambda(\bbP,r)$ be the
following game of two players, {\em Complete} and  {\em Incomplete}:   
\begin{quotation}
\noindent the game lasts $\lambda$ moves and during a play the players
construct a sequence $\langle (p_i,q_i): i<\lambda\rangle$ of pairs of
conditions from $\bbP$ in such a way that $(\forall j<i<\lambda)(r\leq p_j
\leq q_j\leq p_i)$ and at the stage $i<\lambda$ of the game, first
Incomplete chooses $p_i$ and then Complete chooses $q_i$.  
\end{quotation}
Complete wins if and only if for every $i<\lambda$ there are legal moves for
both players. 
\item We say that the forcing notion $\bbP$ is {\em strategically
$({<}\lambda)$--complete\/} if Complete has a winning strategy in the game 
$\Game_0^\lambda(\bbP,r)$ for each condition $r\in\bbP$. 
\item Let $N\prec (\cH(\chi),\in,<^*_\chi)$ be a model such that
${}^{<\lambda} N\subseteq N$, $|N|=\lambda$ and $\bbP\in N$. We say that a
condition $p\in\bbP$ is {\em $(N,\bbP)$--generic in the standard
sense\/} (or just: {\em $(N,\bbP)$--generic\/}) if for every
$\bbP$--name $\name{\tau}\in N$ for an ordinal we have $p\forces$``
$\name{\tau}\in N$ ''. 
\item $\bbP$ is {\em $\lambda$--proper in the standard sense\/} (or just:
{\em $\lambda$--proper\/}) if there is $x\in \cH(\chi)$  such that for
every model $N\prec (\cH(\chi),\in,<^*_\chi)$ satisfying  
\[{}^{<\lambda} N\subseteq N,\quad |N|=\lambda\quad\mbox{ and }\quad\bbP,x
  \in N, \]
and every condition $q\in N\cap\bbP$ there is an $(N,\bbP)$--generic
condition $p\in\bbP$ stronger than $q$.
\end{enumerate}
\end{definition}

\begin{proposition}
[{\cite[Prop. A.1.4]{RoSh:777}}]
\label{obsA.4}
Suppose that $\bbP$ is a $({<}\lambda)$--strategically complete (atomless)
forcing notion, $\alpha^*<\lambda$ and $q_\alpha\in\bbP$ (for $\alpha<
\alpha^*$). Then there are conditions $p_\alpha\in\bbP$ (for $\alpha<
\alpha^*$) such that $q_\alpha\leq p_\alpha$ and for distinct
$\alpha,\alpha'<\alpha^*$ the conditions $p_\alpha, p_{\alpha'}$ are
incompatible.
\end{proposition}

\begin{proposition}
[{\cite[Prop. A.1.6]{RoSh:777}}]
\label{pA.6}
Suppose $\bar{\bbQ}=\langle\bbP_i,\name{\bbQ}_i: i<\gamma\rangle$ is a
$\lambda$--support iteration and, for each $i<\gamma$, 
\[\forces_{\bbP_i}\mbox{`` $\name{\bbQ}_i$ is strategically
$({<}\lambda)$--complete ''.}\]
Then, for each $\vare\leq\gamma$ and $r\in\bbP_\vare$, there
is a winning strategy $\st(\vare,r)$ of Complete in the game
$\Game_0^\lambda(\bbP_\vare,r)$ such that, whenever $\vare_0<\vare_1\leq
\gamma$ and $r\in\bbP_{\vare_1}$, we have:
\begin{enumerate}
\item[(i)]  if $\langle (p_i,q_i):i<\lambda\rangle$ is a play of
$\Game_0^\lambda(\bbP_{\vare_0},r\rest\vare_0)$ in which Complete follows
the strategy $\st(\vare_0,r\rest\vare_0)$, then $\langle (p_i\conc r\rest
[\vare_0,\vare_1),q_i\conc r\rest [\vare_0,\vare_1)):i<\lambda\rangle$ is a
play of $\Game_0^\lambda(\bbP_{\vare_1},r)$ in which Complete uses
$\st(\vare_1,r)$;  
\item[(ii)] if $\langle (p_i,q_i):i<\lambda\rangle$ is a play of
$\Game_0^\lambda(\bbP_{\vare_1},r)$ in which Complete plays according to the
strategy $\st(\vare_1,r)$, then $\langle (p_i\rest\vare_0,q_i\rest\vare_0):i
<\lambda\rangle$ is a play of $\Game_0^\lambda(\bbP_{\vare_0},r\rest
\vare_0)$ in which Complete uses $\st(\vare_0,r\rest\vare_0)$; 
\item[(iii)] is $\vare_1$ is limit and a sequence $\langle (p_i,q_i):
i<\lambda\rangle\subseteq\bbP_{\vare_1}$ is such that for each
$\xi<\vare_1$, $\langle (p_i\rest\xi,q_i\rest\xi):i<\lambda\rangle$ is a
play  of $\Game_0^\lambda(\bbP_\xi,r\rest\xi)$ in which Complete uses the
strategy $\st(\xi,r\rest\xi)$, then $\langle (p_i,q_i):i<\lambda\rangle$ is
a play of $\Game_0^\lambda(\bbP_{\vare_1},r)$ in which Complete plays
according to $\st(\vare_1,r)$; 
\item[(iv)]  if $\langle (p_i,q_i):i<i^*\rangle$ is a partial play of
$\Game_0^\lambda(\bbP_{\vare_1},r)$ in which Complete uses $\st(\vare_1,r)$
and $p'\in\bbP_{\vare_0}$ is stronger than all $p_i\rest\vare_0$ (for
$i<i^*$), then there is $p^*\in \bbP_{\vare_1}$ such that
$p'=p^*\rest\vare_0$ and $p^*\geq p_i$ for $i<i^*$. 
\end{enumerate}
\end{proposition}

\begin{definition}
[{\cite[Def. A.1.7]{RoSh:777}}, see also {\cite[A.3.3, A.3.2]{Sh:587}}]
\label{dA.5}
\begin{enumerate}
\item Let $\alpha,\gamma$ be ordinals, $\emptyset\neq w\subseteq \gamma$.
{\em A standard $(w,\alpha)^\gamma$--tree\/} is a pair $\cT=(T,\rk)$ such
that  
\begin{itemize}
\item $\rk:T\longrightarrow w\cup\{\gamma\}$, 
\item if $t\in T$ and $\rk(t)=\vare$, then $t$ is a sequence $\langle
(t)_\zeta: \zeta\in w\cap\vare\rangle$, where each $(t)_\zeta$ is a
sequence of length $\alpha$,
\item $(T,\vtl)$ is a tree with root $\langle\rangle$ and such that every
chain in $T$ has a $\vartriangleleft$--upper bound it $T$.
\end{itemize}
We will keep the convention that $\cT^x_y$ is $(T^x_y,\rk^x_y)$.
\item Let $\bar{\bbQ}=\langle\bbP_i,\name{\bbQ}_i:i<\gamma\rangle$ be a
$\lambda$--support iteration. {\em A standard tree of conditions in
$\bar{\bbQ}$} is a system $\bar{p}=\langle p_t:t\in T\rangle$ such that 
\begin{itemize}
\item $(T,\rk)$ is a standard $(w,\alpha)^\gamma$--tree for some $w\subseteq
\gamma$ and an ordinal $\alpha$, 
\item $p_t\in\bbP_{\rk(t)}$ for $t\in T$, and
\item if $s,t\in T$, $s\vtl t$, then $p_s=p_t\rest\rk(s)$. 
\end{itemize}
\item Let $\bar{p}^0,\bar{p}^1$ be standard trees of conditions in
  $\bar{\bbQ}$, $\bar{p}^i=\langle p^i_t:t\in T\rangle$. We write
  $\bar{p}^0\leq \bar{p}^1$ whenever for each $t\in T$ we have $p^0_t\leq
  p^1_t$.  
\end{enumerate}
\end{definition}

\begin{proposition}
\label{pA.7}
Assume that $\bar{\bbQ}=\langle\bbP_i,\name{\bbQ}_i:i<\gamma\rangle$ is a
$\lambda$--support iteration such that for all $i<\gamma$ we have 
\[\forces_{\bbP_i}\mbox{`` $\name{\bbQ}_i$ is strategically
$({<}\lambda)$--complete ''.}\]
\begin{enumerate}
\item \cite[Prop. A.1.9]{RoSh:777}\quad Suppose that $\bar{p}=\langle
p_t:t\in T\rangle$ is a standard tree of conditions in $\bar{\bbQ}$,
$|T|<\lambda$, and $\cI\subseteq\bbP_\gamma$ is open dense. Then there is
a standard tree of conditions $\bar{q}=\langle q_t:t\in T\rangle$ such
that $\bar{p}\leq \bar{q}$ and $(\forall t\in T)(\rk(t)=\gamma\
\Rightarrow\ q_t\in\cI)$. 
\item If $\bar{p}=\langle p_t:t\in T\rangle$ is a standard tree of
conditions in $\bar{\bbQ}$ and $|T|<\lambda$, then there is a standard tree 
of conditions $\bar{q}=\langle q_t:t\in T\rangle$ such that $\bar{p}\leq
\bar{q}$ and 
\begin{itemize}
\item if $t_0,t_1\in T$, $\rk(t_0)=\rk(t_1)$, $\xi\in\Dom(t_0)$ and
  $(t_0)_\xi\neq (t_1)_\xi$, $t_0\rest\xi=t_1\rest\xi$ then
\[q_{t_0}\rest\xi\forces_{\bbP_\xi}\mbox{`` the conditions }q_{t_0}(\xi),
  q_{t_1}(\xi) \mbox{ are incompatible in }\name{\bbQ}_\xi\mbox{ ''.}\]
\end{itemize}
\item Suppose that 
\begin{itemize}
\item $w\subseteq\lambda$, $|w|<\lambda$, $1<\mu\leq\lambda$ and
  $T=\bigcup\limits_{\xi\leq\gamma}\prod\limits_{\zeta\in w\cap \xi}\mu$ (so
  $\cT=(T,\rk)$ is a standard $(w,1)^\gamma$--tree),     
\item $\bar{p}=\langle p_t:t\in T\rangle$ is a standard tree of
conditions in $\bar{\bbQ}$, 
\item for $\xi\in w$, $\name{\vare}_\xi$ is a $\bbP_\xi$--name for a
  non-zero ordinal below $\mu$.
\end{itemize}
Then there are a standard $(w,1)^\gamma$--tree $\cT'=(T',\rk')$ and a tree
of conditions $\bar{q}=\langle q_t:t\in T'\rangle$ such that  
\begin{itemize}
\item $T'\subseteq T$ and for every $t\in T'$ such that $\rk'(t)=\xi\in w$
the condition $q_t$ decides the value of $\name{\vare}_\xi$, say $q_t
\forces\name{\vare}_\xi=\vare^t_\xi$, and 
\item $p_t\leq q_t$ for $t\in T'$, and  
\item if $t\in T'$, $\rk(t)=\xi\in w$, then 
\[\big\{\alpha<\mu_\xi: t\cup\{\langle\xi,\alpha\rangle\}\in T'\big\}=
\vare^t_\xi.\]  
\end{itemize}
\end{enumerate}
\end{proposition}

\begin{proof}
(2)\quad Straightforward application of \ref{obsA.4}.
\medskip

\noindent (3)\quad  Note that we cannot apply the first part directly, as
the tree $T$ may be of size $\lambda$. So we will proceed inductively
constructing initial levels of $T'$ of size $<\lambda$ and applying (1) to
them.  

For $\vare\leq\gamma$ and $r\in\bbP_\vare$ let $\st(\vare,r)$ be the winning
strategy of Complete in $\Game^\lambda_0(\bbP_\vare,r)$ given by
\ref{pA.6} (so these strategies have the coherence properties listed there).
Let $\langle\xi_\beta:\beta\leq\beta^*\rangle$ be the increasing enumeration
of $w\cup\{\gamma\}$, $\beta^*<\lambda$. By induction on $\beta\leq\beta^*$
we will pick $\cT_\beta,\bar{q}^\beta,\bar{r}^\beta$ and $\bar{\vare}^\beta$
such that    
\begin{enumerate}
\item[(a)] $\cT_\beta=(T_\beta,\rk_\beta)$ is a standard $(w\cap
\xi_\beta)^\gamma$--tree, $|T_\beta|<\lambda$, and $\bar{q}^\beta=\langle
q^\beta_t:t\in T_\beta\rangle$, $\bar{r}^\beta=\langle r^\beta_t:t\in
T_\beta\rangle$ are tree of conditions, $\bar{q}^\beta\leq\bar{r}^\beta$;   
\item[(b)] if $\beta_0<\beta_1\leq\beta^*$, then $T_{\beta_0}=\{t\rest
\xi_{\beta_0}:t\in T_{\beta_1}\}$ and $r^{\beta_0}_{t\rest\xi_{\beta_0}}\leq 
q^{\beta_1}_t\rest\xi_{\beta_0}$ for $t\in T_{\beta_1}$; 
\item[(c)] if $\beta<\beta^*$, $t\in T_\beta$ and $\rk_\beta(t)=\gamma$ (so
$\rk(t)=\xi_\beta$), then 
\[\langle \big(q^\alpha_{t\rest\xi_\alpha}\!\conc p_t\rest[\xi_\alpha,
\xi_\beta),r^\alpha_{t\rest\xi_\alpha}\!\conc p_t\rest[\xi_\alpha,\xi_\beta) 
\big):\alpha<\beta\rangle^{\textstyle\frown}\langle \big(q^\alpha_t,
r^\alpha_t\big):\beta\leq\alpha<\beta^*\rangle\]
is a partial play of $\Game^\lambda_0(\bbP_{\xi_\beta},p_t)$ in which
Complete uses her winning strategy $\st(\xi_\beta,p_t)$;     
\item[(d)] $\bar{\vare}^\beta=\langle\vare^\beta_t:t\in T_\beta,\
  \rk_\beta(t)=\gamma\rangle\subseteq\lambda$; 
\item[(e)] if $\beta<\beta^*$, $t\in T_\beta$ and $\rk_\beta(t)=\gamma$ (so
$\rk(t)=\xi_\beta$), then $p_t\leq q^\beta_t\in\bbP_{\xi_\beta}$ and 
$q^\beta_t\forces_{\bbP_{\xi_\beta}}\name{\vare}_{\xi_\beta}=
\vare^\beta_t$;  
\item[(f)] if $\beta<\beta^*$, $t\in T_\beta$ and $\rk_\beta(t)=\gamma$,
then $\big\{\alpha<\lambda:t\cup\{\langle\xi_\beta,\alpha\rangle\big\}\in 
  T_{\beta+1}\}=\vare_t^\beta$. 
\end{enumerate}

We let $T_0=\{\langle\rangle\}$ and we choose $q^0_{\langle\rangle}\in
\bbP_{\xi_0}$ and $\vare^0_{\langle\rangle}$ so that $p_{\langle\rangle}\leq
q^0_{\langle\rangle}$ and $q^0_{\langle\rangle}\forces_{\bbP_{\xi_0}}
\name{\vare}_{\xi_0}=\vare^0_{\langle\rangle}$. Then we let $r^0_{\langle
\rangle}$ be the answer given by $\st(\xi_0,p_{\langle\rangle})$. Now
suppose that we have defined $\cT_\alpha,\bar{q}^\alpha,\bar{r}^\alpha$ and
$\bar{\vare}^\alpha$ for $\alpha<\beta\leq\beta^*$. 

If $\beta$ is a limit ordinal then the demands (a) and (b) uniquely define
the standard tree $\cT_\beta$. It follows from the choice of $\st(\vare,r)$
(see clause \ref{pA.6}(iii)) and demand (c) at previous stages that 
\begin{enumerate}
\item[$(\oplus)_\beta$] if $t\in T_\beta$, $\rk_\beta(t)=\gamma$ (so
  $\rk(t)=\xi_\beta$), then the sequence 
\[\big\langle \big(q^\alpha_{t\rest\xi_\alpha}\conc p_t\rest[\xi_\alpha,
\xi_\beta),r^\alpha_{t\rest\xi_\alpha}\conc p_t\rest[\xi_\alpha,\xi_\beta)
\big):\alpha<\beta\big\rangle\]
is a partial play of $\Game^\lambda_0(\bbP_{\xi_\beta},p_t)$ in which
Complete uses her winning strategy $\st(\xi_\beta,p_t)$.
\end{enumerate}
For $t\in T_\beta$ we define a condition $q_t\in \bbP_{\xi_\beta}$ as
follows: 
\begin{itemize}
\item $\Dom(q_t)=\bigcup\limits_{\alpha<\beta}\Dom(r^\alpha_{t\rest
\xi_\alpha})\cup\Dom(p_t)\subseteq\rk(t)$,
\item if $\zeta\in\Dom(q_t)$, then $q_t(\zeta)$ is the $<^*_\chi$--first
$\bbP_\zeta$--name for an element of $\name{\bbQ}_\zeta$ such that
\[\begin{array}{r}
q_t\rest\zeta\forces_{\bbP_\zeta}\mbox{`` if the set }\{r^\alpha_{t\rest
\xi_\alpha}(\zeta):\zeta<\xi_\alpha\ \&\ \alpha<\beta\}\cup\{p_t(\zeta)\}
\mbox{ has an upper bound, }\\
\mbox{\ \ then $q_t(\zeta)$ is such an upper bound ''.}
  \end{array}\]
\end{itemize}
It follows from $(\oplus)_\beta$ (and \ref{pA.6}(iv)) that $p_t\leq q_t$ and
$r^\alpha_{t\rest\xi_\alpha}\leq q_t\rest\xi_{\alpha+1}$ for
$\alpha<\beta$. Now, by ``the $<^*_\chi$--first'', clearly $\bar{q}=\langle
q_t:t\in T_\beta\rangle$ is a tree of conditions. Applying \ref{pA.7}(1) we
may choose a tree of conditions $\bar{q}^\beta=\langle q^\beta_t:t\in
T_\beta\rangle$ such that $\bar{q}\leq\bar{q}^\beta$ and 
\begin{itemize}
\item if $\beta<\beta^*$, $t\in T_\beta$ and $\rk_\beta(t)=\gamma$, then the
condition $q^\beta_t$ decides the value of $\name{\vare}_{\xi_\beta}$ (and
let $q^\beta_t\forces\name{\vare}_{\xi_\beta}=\vare^\beta_t$) and $q^\beta_t
\in \bbP_{\xi_\beta}$.
\end{itemize}
Then, for $t\in T_\beta$, we let $r^\beta_t$ be the answer given to Complete
by $\st(\rk(t),p_t)$ in the appropriate partial play of $\Game^\lambda_0(
\bbP_{\rk(t)},p_t)$, where at stage $\beta$ Incomplete put $q^\beta_t$ (see
(c), $(\oplus)_\beta$). It follows from \ref{pA.6}(ii) that $\bar{r}^\beta=
\langle r^\beta_t:t\in T_\beta\rangle$ is a tree of conditions. Plainly,
$\cT_\beta,\bar{q}^\beta,\bar{r}^\beta$ and $\bar{\vare}^\beta$ satisfy all
relevant (restrictions of the) demands (a)--(f). 

Now suppose that $\beta$ is a successor ordinal, say $\beta=\beta_0+1$. Let 
\[T_\beta=T_{\beta_0}\cup\big\{t\cup\{\langle\xi_{\beta_0},\vare\rangle\}:
t\in T_{\beta_0}\ \&\ \rk_{\beta_0}(t)=\gamma\ \&\ \vare<\vare^{\beta_0}_t
\big\}\]
and for $t\in T_\beta$ define $q_t$ as follows:
\begin{itemize}
\item if $t\in T_{\beta_0}$, then $q_t=r^{\beta_0}_t$,
\item if $t\in T_\beta\setminus T_{\beta_0}$, then $q_t= r^{\beta_0}_{t\rest
\xi_{\beta_0}}\conc p_t\rest [\xi_{\beta_0},\xi_\beta)$.
\end{itemize}
Then $\bar{q}=\langle q_t:t\in T_\beta\rangle$ is a tree of conditions,
$r^{\beta_0}_t\leq q_t$ for $t\in T_{\beta_0}$. It follows from
\ref{pA.7}(1) that we may choose a tree of conditions $\bar{q}^\beta=\langle
q^\beta_t:t\in T_\beta\rangle$ such that $\bar{q}\leq\bar{q}^\beta$ and 
\begin{itemize}
\item if $\beta<\beta^*$, $t\in T_\beta$ and $\rk_\beta(t)=\gamma$, then
the condition  $q^\beta_t$ decides $\name{\vare}_{\xi_\beta}$ and, say,
$q^\beta_t\forces\name{\vare}_{\xi_\beta}=\vare^\beta_t$.
\end{itemize}
Next, like in the limit case, $\bar{r}^\beta=\langle r^\beta_t:t\in
T_\beta\rangle$ is obtained by applying the strategies $\st(\rk(t),p_t)$
suitably. Easily, $\cT_\beta,\bar{q}^\beta,\bar{r}^\beta$ and
$\bar{\vare}^\beta$ satisfy the demands (a)--(f). 

After the inductive construction is carried out look at $T_{\beta^*}$,
$\bar{q}^{\beta^*}$ and $\langle\bar{\vare}^\beta:\beta<\beta^*\rangle$.   
\end{proof}

\section{ABC of reasonable completeness}

\begin{remark}
Note that if $\bbQ$ is strategically $({<}\lambda)$--complete and $\cU$ is a
normal filter on $\lambda$, then the normal filter generated by $\cU$ in
$\bV^{\bbQ}$ is proper. Abusing notation, we may denote the normal filter
generated by $\cU$ in $\bV^\bbQ$ also by $\cU$ or by $\cU^\bbQ$. Thus if
$\name{A}$ is a $\bbQ$--name for a subset of $\lambda$, then $p\forces_\bbQ
\name{A}\in\cU^\bbQ$ if and only if for some $\bbQ$--names $\name{A}_\alpha$
for elements of $\cU^\bV$ we have that $p\forces_\bbQ
\mathop{\triangle}\limits_{\alpha<\lambda}\name{A}_\alpha\subseteq\name{A}$. 
\end{remark}

\begin{definition}
\label{p.1A}
Let $\bbQ$ be a strategically $({<}\lambda)$--complete forcing notion.
\begin{enumerate}
\item For a condition $p\in\bbQ$ we define a game $\Agame(p,\bbQ)$ between
two players, Generic and Antigeneric, as follows. A play of $\Agame(p,\bbQ)$
lasts $\lambda$ steps and during a play a sequence     
\[\Big\langle I_\alpha,\langle p^\alpha_t,q^\alpha_t:t\in I_\alpha\rangle:
\alpha<\lambda\Big\rangle\]
is constructed. Suppose that the players have arrived to a stage $\alpha<
\lambda$ of the game. Now, 
\begin{enumerate}
\item[$(\aleph)_\alpha$]  first Generic chooses a non-empty set $I_\alpha$
of cardinality $<\mu_\alpha$ and a system $\langle p^\alpha_t:t\in I_\alpha
\rangle$ of conditions from $\bbQ$,
\item[$(\beth)_\alpha$]  then Antigeneric answers by picking a system 
$\langle q^\alpha_t:t\in I_\alpha\rangle$ of conditions from $\bbQ$ such that 
$(\forall t\in I_\alpha)(p^\alpha_t\leq q^\alpha_t)$. 
\end{enumerate}
At the end, Generic wins the play 
\[\Big\langle I_\alpha,\langle p^\alpha_t,q^\alpha_t:t\in I_\alpha\rangle:
\alpha<\lambda\Big\rangle\]  
of $\Agame(p,\bbQ)$ \quad if and only if 
\begin{enumerate}
\item[$(\circledast)^{\rm rc}_{\rm A}$] there is a condition $p^*\in\bbQ$
stronger than $p$ and such that\footnote{equivalently, 
for every $\alpha<\lambda$ the set $\big\{q^\alpha_t:t\in
I_\alpha\big\}$ is pre-dense above $p^*$}
\[p^*\forces_{\bbQ}\mbox{`` }\big\{\alpha<\lambda:\big(\exists t\in I_\alpha
\big)\big(q^\alpha_t\in\Gamma_{\bbQ}\big)\big\}=\lambda\mbox{ ''}.\]
\end{enumerate}
\item Games $\tegame(p,\bbQ),\Cgame(p,\bbQ)$ are defined similarly, except
that the winning criterion $(\circledast)^{\rm rc}_{\rm A}$ is replaced by 
\begin{enumerate}
\item[$(\circledast)^{\rm rc}_{\rm B}$] there is a condition $p^*\in\bbQ$
stronger than $p$ and such that  
\[p^*\forces_{\bbQ}\mbox{`` }\big\{\alpha<\lambda:\big(\exists t\in I_\alpha
\big)\big(q^\alpha_t\in\Gamma_{\bbQ}\big)\big\}\in \cU^\bbQ\mbox{ ''},\]
\item[$(\circledast)^{\rm rc}_{\rm C}$] there is a condition $p^*\in\bbQ$
stronger than $p$ and such that  
\[p^*\forces_{\bbQ}\mbox{`` }\big\{\alpha<\lambda:\big(\exists t\in I_\alpha
\big)\big(q^\alpha_t\in\Gamma_{\bbQ}\big)\big\}\in \big(\cU^\bbQ\big)^+
\mbox{ ''},\]
\end{enumerate}
respectively.
\item For a condition $p\in\bbQ$ we define a game $\bgame(p,\bbQ)$ between
Generic and Antigeneric as follows. A play of $\bgame(p,\bbQ)$ lasts
$\lambda$ steps and during a play a sequence      
\[\Big\langle \zeta_\alpha,\langle p^\alpha_\xi,q^\alpha_\xi:\xi<
\zeta_\alpha\rangle:\alpha<\lambda\Big\rangle\]
is constructed. Suppose that the players have arrived to a stage $\alpha<
\lambda$ of the game. Now, Generic chooses a non-zero ordinal
$\zeta_\alpha<\mu_\alpha$ and then the two players play a subgame of length  
$\zeta_\alpha$ alternatively choosing successive terms of a sequence 
$\langle p^\alpha_\xi,q^\alpha_\xi:\xi<\zeta_\alpha\rangle$. At a stage
$\xi<\zeta_\alpha$ of the subgame, first Generic picks a condition
$p^\alpha_\xi\in\bbQ$ and then Antigeneric answers with a condition
$q^\alpha_\xi$ stronger than $p^\alpha_\xi$.

At the end, Generic wins the play 
\[\Big\langle \zeta_\alpha,\langle p^\alpha_\xi,q^\alpha_\xi:\xi<
\zeta_\alpha\rangle:\alpha<\lambda\Big\rangle\]  
of $\bgame(p,\bbQ)$\quad if and only if 
\begin{enumerate}
\item[$(\circledast)^{\rm rc}_{\bf b}$] there is a condition $p^*\in\bbQ$ 
stronger than $p$ and such that  
\[p^*\forces_{\bbQ}\mbox{`` }\big\{\alpha<\lambda:\big(\exists
\xi<\zeta_\alpha\big)\big(q^\alpha_\xi\in\Gamma_{\bbQ}\big)\big\}\in
\cU^\bbQ \mbox{ ''}.\] 
\end{enumerate}
\item Games $\agame(p,\bbQ)$ and $\cgame(p,\bbQ)$ are defined similarly
  except that the winning criterion $(\circledast)^{\rm rc}_{\bf b}$ is
  changed so that ``$\in \cU^\bbQ$'' is replaced by ``$=\lambda$'' or ``$\in
  \big(\cU^\bbQ\big)^+$'', respectively. 
\item We say that a forcing notion $\bbQ$ is {\em reasonably A--bounding 
over $\bar{\mu}$\/} if    
\begin{enumerate}
\item[(a)] $\bbQ$ is strategically $({<}\lambda)$--complete,  and 
\item[(b)] for any $p\in\bbQ$, Generic has a winning strategy in the game
$\Agame(p,\bbQ)$. 
\end{enumerate}
In an analogous manner we define when the forcing notion $\bbQ$ is {\em
reasonably X--bounding over $\cU,\bar{\mu}$\/} (for ${\rm X}\in\{{\rm B},
{\rm C}, {\bf a}, {\bf b}, {\bf c}\}$) --- just using the game $\Game^{\rm 
rcX}_{\cU,\bar{\mu}}(p,\bbQ)$ appropriately. 

If $\mu_\alpha=\lambda$ for each $\alpha<\lambda$, then we may omit
$\bar{\mu}$ and say {\em reasonably B--bounding over $\cU$\/} etc. If $\cU$ 
is the filter generated by club subsets of $\lambda$, we may omit it as
well. 
\item Let $\st$ be a strategy for Generic in the game $\tegame(p,\bbQ)$. We
will say that a sequence $\Big\langle I_\alpha,\langle p^\alpha_t,
q^\alpha_t:t\in I_\alpha\rangle:\delta<\alpha<\lambda\Big\rangle$ is {\em 
a $\delta$--delayed play according to $\st$} if it has an extension $\Big
\langle I_\alpha,\langle p^\alpha_t,q^\alpha_t:t\in I_\alpha\rangle:\alpha<
\lambda\Big\rangle$ which is a play agreeing with $\st$ and such that
$p^\alpha_t=q^\alpha_t$ for $\alpha\leq\delta$, $t\in I_\alpha$. 
\end{enumerate}
\end{definition}

\begin{remark}
If $\st$ is a winning strategy for Generic in the game $\tegame(p,\bbQ)$,
and $\bar{\sigma}=\Big\langle I_\alpha,\langle p^\alpha_t, q^\alpha_t:t
\in I_\alpha\rangle:\delta\leq\alpha<\lambda\Big\rangle$ is a
$\delta$--delayed play according to $\st$, then $\bar{\sigma}$ satisfies the
condition $(\circledast)^{\rm rc}_{\rm B}$.  
\end{remark}

\begin{observation}
For $\cU,\bar{\mu}$ as in \ref{incon}, $X\in\{A,B,C,{\bf a},{\bf b},{\bf
  c}\}$ and a forcing notion $\bbQ$, let $\Phi(\bbQ,X,\cU,\bar{\mu})$ be the
  statement  
\begin{center}
``$\bbQ$ is reasonably $X$--bounding over $\cU,\bar{\mu}$''. 
\end{center}
Then the following implications hold

\[\begin{array}{ccccccc}
\Phi(\bbQ,A,\bar{\mu})& \Rightarrow & \Phi(\bbQ,B,\cU,\bar{\mu}) &
\Rightarrow & \Phi(\bbQ,C,\cU,\bar{\mu}) & & \\
\Downarrow & & \Downarrow & & \Downarrow & & \\
\Phi(\bbQ,{\bf a},\bar{\mu})& \Rightarrow & \Phi(\bbQ,{\bf b},\cU,
\bar{\mu}) & \Rightarrow & \Phi(\bbQ,{\bf c},\cU,\bar{\mu})& \Rightarrow &   
\bbQ\mbox{ is $\lambda$--proper.}
  \end{array}\]
\end{observation}

\begin{theorem}
\label{verB}
Assume that $\lambda,\cU,\bar{\mu}$ are as in \ref{incon} and $\bar{\bbQ}= 
\langle\bbP_\xi,\name{\bbQ}_\xi:\xi<\gamma\rangle$ is a $\lambda$--support
iteration such that for every $\xi<\gamma$,  
\[\forces_{\bbP_\xi}\mbox{`` $\name{\bbQ}_\xi$ is reasonably B--bounding
  over $\cU,\bar{\mu}$ ''.}\] 
Then $\bbP_{\gamma}=\lim(\bar{\bbQ})$ is reasonably ${\bf b}$--bounding over  
$\cU,\bar{\mu}$ (and so also $\lambda$--proper). 
\end{theorem}

\begin{proof}
For each $\xi<\gamma$ pick a $\bbP_\xi$--name $\name{\st}^0_\xi\in N$ such
that  
\[\begin{array}{r}
\forces_{\bbP_\xi}\mbox{`` }\name{\st}^0_\xi\mbox{ is a winning strategy
  for Complete in }\Game_0^\lambda\big(\name{\bbQ}_\xi,
  \name{\emptyset}_{\name{\bbQ}_\xi}\big)\mbox{ such that }\ \\
\mbox{ if Incomplete plays $\name{\emptyset}_{\name{\bbQ}_\xi}$ then
  Complete answers with $\name{\emptyset}_{\name{\bbQ}_\xi}$ as well ''.}
  \end{array}\] 
Also, for $\xi\leq\gamma$ and $r\in\bbP_\xi$, let $\st(\xi,r)$ be a winning
strategy of Complete in $\Game_0^\lambda(\bbP_\xi,r)$ with the coherence
properties given in \ref{pA.6}. 

We are going to describe a strategy $\st$ for Generic in the game
$\bgame(p,\bbP_\gamma)$. In the course of the play, at a stage
$\delta<\lambda$, Generic will be instructed to construct aside
\begin{enumerate}
\item[$(\otimes)_\delta$] \qquad $\cT_\delta,\bar{p}^\delta_*,
\bar{q}^\delta_*,r^-_\delta, r_\delta,w_\delta$, $\langle
\name{\vare}_{\delta,\xi},\name{\bar{p}}_{\delta,\xi},
\name{\bar{q}}_{\delta,\xi}:\xi\in w_\delta\rangle$, and $\name{\st}_\xi$
for $\xi\in w_{\delta+1}\setminus w_\delta$.  
\end{enumerate}
These objects will be chosen so that if 
\[\big\langle\zeta_\delta,\langle p^\delta_\zeta,q^\delta_\zeta:\zeta<
\zeta_\delta\rangle:\delta<\lambda\big\rangle\] 
is a play of $\bgame(p,\bbP_\gamma)$ in which Generic follows $\st$, and the
objects constructed at stage $\delta<\lambda$ are listed in
$(\otimes)_\delta$, then the following conditions are satisfied (for each
$\delta<\lambda$).  
\begin{enumerate}
\item[$(*)_1$] $r^-_\delta,r_\delta\in \bbP_\gamma$, $r_0(0)=p(0)$,
$w_\delta\subseteq\gamma$, $|w_\delta|=|\delta|+1$, $\bigcup\limits_{\alpha< 
\lambda}\Dom(r_\alpha)=\bigcup\limits_{\alpha<\lambda} w_\alpha$,
$w_0=\{0\}$, $w_\delta\subseteq w_{\delta+1}$ and if $\delta$ is limit then
$w_\delta=\bigcup\limits_{\alpha<\delta} w_\alpha$.
\item[$(*)_2$] For each $\alpha<\delta<\lambda$ we have $(\forall\xi\in
w_{\alpha+1})(r_\alpha(\xi)=r_\delta(\xi))$ and $p\leq r_\alpha^-\leq
r_\alpha\leq r^-_\delta\leq r_\delta$. 
\item[$(*)_3$] If $\xi\in\gamma\setminus w_\delta$, then 
\[\begin{array}{ll}
r_\delta\rest\xi\forces&\mbox{`` the sequence }\langle r^-_\alpha(\xi),
r_\alpha(\xi):\alpha\leq\delta\rangle\mbox{ is a legal partial play of }\\
&\quad\Game_0^\lambda\big(\name{\bbQ}_\xi,\name{\emptyset}_{\name{\bbQ}_\xi} 
\big)\mbox{ in which Complete follows }\name{\st}^0_\xi\mbox{ ''} 
\end{array}\]
and if $\xi\in w_{\delta+1}\setminus w_\delta$, then $\name{\st}_\xi$ is a
$\bbP_\xi$--name  for a winning strategy of Generic in
$\tegame(r_\delta(\xi),\name{\bbQ}_\xi)$ such that if $\langle p^\alpha_t:
t\in I_\alpha\rangle$ is given by that strategy to Generic at stage
$\alpha$, then $I_\alpha$ is an ordinal below $\mu_\alpha$. (And $\st_0$ is
a suitable winning strategy of Generic in $\tegame(p(0),\bbQ_0)$.)   
\item[$(*)_4$] $\cT_\delta=(T_\delta,\rk_\delta)$ is a standard $(w_\delta,
1)^\gamma$--tree, $|T_\delta|<\mu_\delta$.   
\item[$(*)_5$] $\bar{p}^\delta_*=\langle p^\delta_{*,t}:t\in T_\delta
\rangle$ and $\bar{q}^\delta_*=\langle q^\delta_{*,t}:t\in T_\delta\rangle$
are standard trees of conditions, $\bar{p}^\delta_*\leq\bar{q}^\delta_*$.  
\item[$(*)_6$] For $t\in T_\delta$ we have $\big(\bigcup\limits_{\alpha<
\delta}\Dom(r_\alpha)\cup w_\delta\big)\cap\rk_\delta(t)\subseteq
  \Dom(p^\delta_{*,t})$ and for each $\xi\in \Dom(p^\delta_{*,t})\setminus
  w_\delta$:
\[\begin{array}{ll}
p^\delta_{*,t}\rest\xi\forces&\mbox{`` if the set }\{r_\alpha(\xi):\alpha<
\delta\}\mbox{ has an upper bound in }\name{\bbQ}_\xi,\\ 
&\mbox{\quad then $p^\delta_{*,t}(\xi)$ is such an upper bound ''.}
  \end{array}\]
\item[$(*)_7$] $\zeta_\delta=|\{t\in T_\delta:\rk_\delta(t)=\gamma\}|$ and
  for some enumeration $\{t\in T_\delta:\rk_\delta(t)=\gamma\}=\{t_\zeta:
  \zeta<\zeta_\delta\}$, for each $\zeta<\zeta_\delta$ we have 
\[p^\delta_{*,t_\zeta}\leq p^\delta_\zeta\leq q^\delta_\zeta\leq
  q^\delta_{*,t_\zeta}.\] 
\item[$(*)_8$] If $\xi\in w_\delta$, then $\name{\vare}_{\delta,\xi}$ is a
  $\bbP_\xi$--name for an ordinal below $\mu_\delta$,
  $\name{\bar{p}}_{\delta,\xi},\name{\bar{q}}_{\delta,\xi}$ are
  $\bbP_\xi$--names for sequences of conditions in $\name{\bbQ}_\xi$ of
  length $\name{\vare}_{\delta,\xi}$.   
\item[$(*)_9$] If $\xi\in w_{\beta+1}\setminus w_\beta$, $\beta<\lambda$,
  then 
\[\begin{array}{r}
\forces_{\bbP_\xi}\mbox{`` }\langle\name{\vare}_{\alpha,\xi},
\name{\bar{p}}_{\alpha,\xi},\name{\bar{q}}_{\alpha,\xi}:\beta<\alpha<
\lambda\rangle\mbox{ is a delayed play of }\tegame(r_\beta(\xi),
\name{\bbQ}_\xi)\\
\mbox{ in which Generic uses $\name{\st}_\xi$ ''.}
  \end{array}\]
\item[$(*)_{10}$] If $t\in T_\delta$, $\rk_\delta(t)=\xi<\gamma$, then
the condition $p^\delta_{*,t}$ decides the value of $\name{\vare}_{\delta,
\xi}$, say $p^\delta_{*,t}\forces$``$\name{\vare}_{\delta,\xi}=
\vare^t_{\delta,\xi}$'', and $\{(s)_\xi:t\vartriangleleft s\in T_\delta\}=
\vare^t_{\delta,\xi}$ and  
\[q^\delta_{*,t}\forces_{\bbP_\xi}\mbox{`` } \name{\bar{p}}_{\delta,\xi}
(\vare)\leq p^\delta_{*,t\conc\langle\vare\rangle}(\xi)\mbox{ for }\vare<
\vare^t_{\delta,\xi}\mbox{ and }\name{\bar{q}}_{\delta,\xi}=\langle
q^\delta_{*,s}(\xi):t\vartriangleleft s\in T_\delta\rangle\mbox{ ''.}\]
\item[$(*)_{11}$] If $t_0,t_1\in T_\delta$, $\rk_\delta(t_0)=
\rk_\delta(t_1)$ and $\xi\in w_\delta\cap\rk_\delta(t_0)$,
$t_0\rest\xi=t_1\rest\xi$ but $\big(t_0\big)_\xi\neq \big(t_1\big)_\xi$,
then 
\[q^\delta_{*,t_0\rest\xi}\forces_{\bbP_\xi}\mbox{`` the conditions  
$q^\delta_{*,t_0}(\xi),q^\delta_{*,t_1}(\xi)$ are incompatible ''.}\] 
\item[$(*)_{12}$] $\Dom(r_\delta)=\bigcup\limits_{t\in T_\delta}
\Dom(q^\delta_{*,t})\cup\Dom(p)$ and if $t\in T_\delta$,
$\xi\in\Dom(r_\delta)\cap \rk_\delta(t)\setminus w_\delta$, and
$q^\delta_{*,t}\rest\xi\leq q\in\bbP_\xi$, $r_\delta\rest\xi\leq q$, then   
\[\begin{array}{ll}
q\forces_{\bbP_\xi}&\mbox{`` if the set }\{r_\alpha(\xi):\alpha<\delta\}
\cup\{q^\delta_{*,t}(\xi), p(\xi)\}\mbox{ has an upper bound in }  
\name{\bbQ}_\xi,\\ 
&\mbox{\quad then $r_\delta(\xi)$ is such an upper bound ''.}
  \end{array}\]
\end{enumerate}

To describe the instructions given by $\st$ at stage $\delta<\lambda$ of a
play of $\bgame(p,\bbP_\gamma)$ let us assume that 
\[\big\langle\zeta_\alpha,\langle p^\alpha_\zeta,q^\alpha_\zeta:\zeta<
\zeta_\alpha\rangle:\alpha<\delta\big\rangle\] 
is the result of the play so far and that Generic constructed objects listed
in $(\otimes)_\alpha$ (for $\alpha<\delta$) with properties
$(*)_1$--$(*)_{12}$. 

First, Generic uses her favourite bookkeeping device to determine $w_\delta$
such that the demands in $(*)_1$ are satisfied (and that at the end we will
have $\bigcup\limits_{\alpha<\lambda}\Dom(r_\alpha)=\bigcup\limits_{\alpha<
\lambda}w_\alpha$). Now Generic lets $\cT_\delta'$ be a standard
$(w_\delta,1)^\gamma$--tree such that for each $\xi\in w_\delta\cup
\{\gamma\}$ we have $\{t\in T_\delta':\rk_\delta'(t)=\xi\}=
\prod\limits_{\vare\in w_\delta\cap\xi}\mu_\delta$. Then for $\xi\in
w_\delta$ she chooses $\bbP_\xi$--names $\name{\vare}_{\delta,\xi},
\name{\bar{p}}_{\delta,\xi}$ such that $\name{\vare}_{\delta,\xi}$ is a name
for an ordinal below $\mu_\delta$ and $\name{\bar{p}}_{\delta,\xi}$ is a
name for a sequence of conditions in $\name{\bbQ}_\xi$ of length
$\name{\vare}_{\delta,\xi}$ and    
\[\begin{array}{ll}
\forces_{\bbP_\xi}&\mbox{`` }\name{\vare}_{\delta,\xi},
\name{\bar{p}}_{\delta,\xi} \mbox{ is the answer to the delayed play }\\ 
&\quad\langle\name{\vare}_{\alpha,\xi},\name{\bar{p}}_{\alpha,\xi}, 
\name{\bar{q}}_{\alpha,\xi}:\xi\in w_\alpha\ \&\ \alpha<\delta\rangle  
\mbox{ given to Complete by }\name{\st}_\xi\mbox{ ''}.
  \end{array}\]
She lets $\bar{p}^{\delta,0}_*=\langle p^{\delta,0}_{*,t}:t\in T'_\delta
\rangle$ be a tree of conditions defined so that $\Dom(p^{\delta,0}_{*,t})=
\big(\bigcup\limits_{\alpha<\delta}\Dom(r_\alpha)\cup w_\delta\big)\cap
\rk_\delta'(t)$ and for each $\xi\in \Dom(p^{\delta,0}_{*,t})$ 
\begin{enumerate}
\item[$(*)_{13}$]  $p^{\delta,0}_{*,t}(\xi)$ is the $<^*_\chi$--first
$\bbP_\xi$--name for an element of $\name{\bbQ}_\xi$ such that 
\begin{itemize}
\item if $\xi\in w_\delta$, then 
\[\forces_{\bbP_\xi}\mbox{`` if }(t)_\xi<\name{\vare}_{\delta,\xi}\mbox{
then }p^{\delta,0}_{*,t}(\xi)=\name{\bar{p}}_{\delta,\xi}\big((t)_\xi\big), 
\mbox{ otherwise }p^{\delta,0}_{*,t}(\xi)=\name{\emptyset}_{\name{\bbQ}_\xi}
\mbox{ '',}\]
\item if $\xi\notin w_\delta$, then 
\[\begin{array}{ll}
\forces_{\bbP_\xi}&\mbox{`` if the set }\{r_\alpha(\xi):\alpha<\delta\}
\mbox{ has an upper bound in }\name{\bbQ}_\xi,\\ 
&\mbox{\quad then $p^{\delta,0}_{*,t}(\xi)$ is such an upper bound ''.}
  \end{array}\]
\end{itemize}
\end{enumerate}
Now Generic uses \ref{pA.7}(3) and then \ref{pA.7}(2) to choose a standard
tree $(w_\delta,1)^\gamma$--tree $\cT_\delta=(T_\delta,\rk_\delta)$ and a
tree of conditions $\bar{p}^\delta_*=\langle p^\delta_{*,t}:t\in T_\delta
\rangle$ such that    
\begin{enumerate}
\item[$(*)_{14}^{\rm a}$] $T_\delta\subseteq T_\delta'$ and for every $t\in
  T_\delta$ such that $\rk_\delta(t)=\xi\in w_\delta$ the condition
  $p^\delta_{*,t}$ decides the value of $\name{\vare}_{\delta,\xi}$, say
  $p^\delta_{*,t}\forces\name{\vare}_{\delta,\xi}=\vare^t_{\delta,\xi}$, and 
\item[$(*)_{14}^{\rm b}$] if $t\in T_\delta$, $\rk_\delta(t)=\xi\in
  w_\delta$, then $\{\alpha<\lambda:t\cup\{\langle\xi,\alpha\rangle\}\in
  T_\delta\}=\vare^t_{\delta,\xi}$, and 
\item[$(*)_{14}^{\rm c}$] $p^{\delta,0}_{*,t}\leq p^\delta_{*,t}$ for all
  $t\in T_\delta$, and if $t_0,t_1\in T_\delta$, $\rk_\delta(t_0)=
  \rk_\delta(t_1)$, $\xi\in\Dom(t_0)$, and $t_0\rest\xi=t_1\rest\xi$ but
  $(t_0)_\xi\neq (t_1)_\xi$, then 
\[p^\delta_{*,t_0\rest\xi}\forces_{\bbP_\xi}\mbox{`` the conditions }
  p^\delta_{*,t_0}(\xi),p^\delta_{*,t_1}(\xi)\mbox{ are incompatible in
  $\name{\bbQ}_\xi$ '',}\]
\end{enumerate}
Thus Generic has written aside $\cT_\delta$, $\bar{p}^\delta_*$, $w_\delta$
and $\langle \name{\vare}_{\delta,\xi},\name{\bar{p}}_{\delta,\xi}:\xi\in
w_\delta\rangle$. (It should be clear that they satisfy the demands in
$(*)_1$, $(*)_4$--$(*)_6$, $(*)_8$ and $(*)_9,(*)_{10}$.) Now she turns to
the play of $\bgame(p,\bbP_\gamma)$ and she puts
\[\zeta_\delta=|\{t\in T_\delta:\rk_\delta(t)=\gamma\}|\]
and she also picks an enumeration $\langle t_\zeta:\zeta<\zeta_\delta
\rangle$ of $\{t\in T_\delta:\rk_\delta(t)=\gamma\}$. The two players start
playing the subgame of level $\delta$ of length $\zeta_\delta$. During the
subgame Generic constructs partial plays $\langle(r^\zeta_i,s^\zeta_i):i
\leq\zeta_\delta\rangle$ of $\Game_0^\lambda(\bbP_\gamma,p^\delta_{*,
t_\zeta})$ (for $\zeta<\zeta_\delta$) in which Complete uses the strategy 
$\st(\gamma,p^\delta_{*,t_\zeta})$ and such that 
\begin{enumerate}
\item[$(*)^{\rm a}_{15}$]  if   $\zeta,\xi<\zeta_\delta$, $t\in T_\delta$,
  $t\vtl t_\zeta$, $t\vtl t_\xi$, $i\leq\zeta_\delta$, then $r^\zeta_i\rest
  \rk_\delta(t)=r^\xi_i\rest\rk_\delta(t)$ and $s^\zeta_i\rest\rk_\delta(t)=
  s^\xi_i\rest \rk_\delta(t)$; 
\item[$(*)^{\rm b}_{15}$]  if $p^\delta_\zeta,q^\delta_\zeta$ are the
  conditions played at stage $\zeta$ of the subgame, then
  $p^\delta_{*,t_\zeta}\leq r^\zeta_i\leq p^\delta_\zeta\leq q^\delta_\zeta
  =r^\zeta_\zeta$ for all $i<\zeta$. 
\end{enumerate}
So suppose that the two players  have arrived at a stage
$\zeta<\zeta_\delta$ of the subgame and $\big\langle\langle (r^\xi_i,
s^\xi_i):i<\zeta\rangle: \xi<\zeta_\delta\big\rangle$ has been
defined. Generic looks at $\langle (r^\zeta_i,s^\zeta_i):i<\zeta\rangle$ --
it is a play of $\Game_0^\lambda(\bbP_\gamma,p^\delta_{*,t_\zeta})$ in which 
Complete uses $\st(\gamma,p^\delta_{*,t_\zeta})$, so we may find a condition
$p^\delta_\zeta\in\bbP_\gamma$ stronger than all $r^\zeta_i,s^\zeta_i$ for
$i<\zeta$ (and $p^\delta_\zeta\geq p^\delta_{*,t_\zeta}$). She plays this
condition as her move at stage $\zeta$ of the subgame and Antigeneric
answers with $q^\delta_\zeta\geq p^\delta_\zeta$. Generic lets
$r^\zeta_\zeta=q^\delta_\zeta$ and she defines $r^\xi_\zeta$ for $\xi<
\zeta_\delta$, $\xi\neq \zeta$, as follows. Let $t\in T_\delta$ be such that
$t\vtl t_\zeta$, $t\vtl t_\xi$ and $\rk_\delta(t)$ is the largest
possible. Generic declares that   
\[\Dom(r^\xi_\zeta)=\big(\Dom(r^\zeta_\zeta)\cap\rk_\delta(t)\big)\cup
\bigcup\limits_{i<\zeta}\Dom(s_i^\xi)\cup\Dom(p^\delta_{*,t_\xi}),\] 
and $r^\xi_\zeta\rest\rk_\delta(t)=r^\zeta_\zeta\rest\rk_\delta(t)$, and for
$\vare\in \Dom(r^\xi_\zeta)\setminus\rk_\delta(t)$ she lets
$r^\xi_\zeta(\vare)$ be the $<^*_\chi$--first $\bbP_\vare$--name for a
member of $\name{\bbQ}_\vare$ such that 
\[r^\xi_\zeta\rest\vare\forces_{\bbP_\vare}\mbox{`` $r^\xi_\zeta(\vare)$ is
  an upper bound to }\{p^\delta_{*,t_\xi}(\vare)\}\cup \{s^\xi_i(\vare):
i<\zeta\}\mbox{ ''}\] 
(remember \ref{pA.6}). Finally,  $s^\xi_\zeta$ (for $\xi<\zeta_\delta$) is
defined as the condition given to Complete by $\st(\gamma,p^\delta_{*,
  t_\zeta})$ in answer to $\langle(r^\xi_i,s^\xi_i):i<\zeta\rangle\conc 
\langle r^\xi_\zeta\rangle$. 

After the subgame is completed and both $p^\delta_\zeta,q^\delta_\zeta$ and
$\big\langle\langle (r^\xi_i,s^\xi_i):i<\zeta_\delta\rangle:\xi<
\zeta_\delta \big\rangle$ has been determined, Generic chooses
$r^0_{\zeta_\delta}$ as any upper bound to $\langle s^0_i:i<\zeta_\delta
\rangle$ and then defines $r^\xi_{\zeta_\delta}$ for
$\xi\in\zeta_\delta\setminus 1$ like $r^\xi_\zeta$ for $\xi\neq \zeta$
above. Also $s^\xi_{\zeta_\delta}$ (for $\xi<\zeta_\delta$) are chosen like
earlier (as results of applying $\st(\gamma,p^\delta_{*,t_\xi})$). Finally,
Generic picks a standard tree of conditions $\bar{q}^\delta_*=\langle
q^\delta_{*,t}:t\in T_\delta\rangle$ such that $(\forall\zeta<\zeta_\delta)
(q^\delta_{*,t_\zeta}= s^\zeta_{\zeta_\delta})$. (Note that $(*)_5$, $(*)_7$
hold.) 

Now Generic defines $r^-_\delta,r_\delta\in\bbP_\gamma$ so that
\[\Dom(r^-_\delta)=\Dom(r_\delta)=\bigcup\limits_{t\in T_\delta}\Dom(
q^\delta_{*,t})\cup \Dom(p)\]
 and 
\begin{enumerate}
\item[$(*)_{16}^{\rm a}$] if $\xi\in\Dom(r^-_\delta)\setminus w_\delta$,
  then:\\
$r^-_\delta(\xi)$ is the $<^*_\chi$--first $\bbP_\xi$--name for an element
  of $\name{\bbQ}_\xi$ such that 
\[\begin{array}{ll}
r^-_\delta\rest\xi\forces_{\bbP_\xi}&\mbox{`` }r^-_\delta(\xi)\mbox{ is an
upper bound of }\{r_\alpha(\xi):\alpha<\delta\}\cup\{p(\xi)\}\mbox{ and }\\  
&\mbox{\quad if }t\in T_\delta,\ \ \rk_\delta(t)>\xi,\mbox{ and }
q^\delta_{*,t}\rest\xi\in\Gamma_{\bbP_\xi}\mbox{ and the set}\\
&\quad\{r_\alpha(\xi):\alpha<\delta\}\cup\{q^\delta_{*,t}(\xi),
  p(\xi)\}\mbox{ has an upper bound in }\name{\bbQ}_\xi,\\ 
&\mbox{\quad then $r^-_\delta(\xi)$ is such an upper bound '',}
  \end{array}\]
and 
$r_\delta(\xi)$ is the $<^*_\chi$--first $\bbP_\xi$--name for an element
  of $\name{\bbQ}_\xi$ such that 
\[\begin{array}{ll}
r_\delta\rest\xi\forces_{\bbP_\xi}&\mbox{`` }r_\delta(\xi)\mbox{ is given to
  Complete by }\name{\st}^0_\xi\mbox{ as the answer to }\ \\ 
&\quad\langle r^-_\alpha(\xi),r_\alpha(\xi):\alpha<\delta\rangle\conc
  \langle r^-_\delta(\xi)\rangle \mbox{  ''} 
  \end{array}\]
\item[$(*)_{16}^{\rm b}$] if $\xi\in w_{\alpha+1}$, $\alpha<\delta$, then
  $r^-_\delta(\xi)=r_\delta(\xi)=r_\alpha(\xi)$. 
\end{enumerate}
(Note that by a straightforward induction on $\xi\in\Dom(r_\delta)$ one
easily applies $(*)_3$ from previous stages to show that
$r_\delta^-,r_\delta$ are well defined and $r_\delta\geq r^-_\delta\geq
r_\alpha,p$ for $\alpha<\delta$. Remember also $(*)_{11}$ and/or $(*)^{\rm
  c}_{14}$.) If $\delta=0$ we also stipulate $r^-_0(0)=r_0(0)=p(0)$.   

Finally, for each $\xi\in w_\delta$, Generic chooses a $\bbP_\xi$--name 
$\name{\bar{q}}_{\delta,\xi}$ for a sequence of conditions in
$\name{\bbQ}_\xi$ of length $\name{\vare}_{\delta,\xi}$ such that  
\[\begin{array}{ll}
\forces_{\bbP_\xi}&\mbox{`` }(\forall\vare<\name{\vare}_{\delta,\xi})(
\name{\bar{p}}_{\delta,\xi}(\vare)\leq\name{\bar{q}}_{\delta,\xi}(\vare))
\mbox{ and }\\
&\mbox{\quad if }t\in T_\delta,\ \rk_\delta(t)>\xi,\mbox{ and }
q^\delta_{*,t}\rest\xi\in\Gamma_{\bbP_\xi}\mbox{ then }
\name{\bar{q}}_{\delta,\xi}\big((t)_\xi\big)=q^\delta_{*,t}(\xi)
\mbox{  ''.}
  \end{array}\]
Generic also picks $w_{\delta+1}$ by the bookkeeping device mentioned at the
beginning and for $\xi\in w_{\delta+1}\setminus w_\delta$ she fixes
$\name{\st}_\xi$  as in $(*)_3$. 
\smallskip

This completes the description of the side objects constructed by Generic at
stage $\delta$. Verification that they satisfy our demands
$(*)_1$--$(*)_{12}$ is straightforward, and thus the description of the
strategy $\st$ is complete.
\medskip

We are going to argue now that $\st$ is a winning strategy for Generic. To
this end suppose that 
\[\big\langle\zeta_\delta,\langle p^\delta_\zeta,q^\delta_\zeta:\zeta<
\zeta_\delta\rangle:\delta<\lambda\big\rangle\] 
is the result of a play of $\bgame(p,\bbP_\gamma)$ in which Generic
followed $\st$ and constructed aside objects listed in $(\otimes)_\delta$
(for $\delta<\lambda$)  so that $(*)_1$--$(*)_{12}$ hold. 

We define a condition $r\in\bbP_\gamma$ as follows. Let $\Dom(r)=
\bigcup\limits_{\delta<\lambda}\Dom(r_\delta)$ and for $\xi\in\Dom(r)$ let
$r(\xi)$ be a $\bbP_\xi$--name for a condition in $\name{\bbQ}_\xi$ such that 
if $\xi\in w_{\alpha+1}\setminus w_\alpha$, $\alpha<\lambda$ (or
$\xi=0=\alpha$), then  
\[\forces_{\bbP_\xi}\mbox{`` }r(\xi)\geq r_{\alpha}(\xi)\mbox{ and }
r(\xi) \forces_{\name{\bbQ}_\xi}\big\{\delta\!<\!\lambda\!:\big(\exists
\vare\!<\!\name{\vare}_{\delta,\xi}\big)\big(\bar{q}_{\delta,\xi}(\vare)\in
\Gamma_{\name{\bbQ}_\xi}\big)\big\}\in (\cU^{\bbP_\xi})^{\name{\bbQ_\xi}}
\mbox{ ''.}\] 
Clearly $r$ is well defined (remember $(*)_9$) and $(\forall\delta<
\lambda)(r_\delta\leq r)$ and $r\geq p$. For each $\xi\in\Dom(r)$ choose a
sequence $\langle\name{A}^\xi_i:i<\lambda\rangle$ of $\bbP_{\xi+1}$--names
for elements of $\cU\cap\bV$ such that   
\begin{enumerate}
\item[$(*)^\xi_{17}$] \quad $\displaystyle r\rest(\xi+1)\forces_{\bbP_{\xi+
1}}\big(\forall\delta\in \mathop{\triangle}\limits_{i<\lambda}
\name{A}^\xi_i\big)\big(\exists\vare<\name{\vare}_{\delta,\xi}\big)\big(
\name{\bar{q}}_{\delta,\xi}(\vare)\in\Gamma_{\name{\bbQ}_\xi}\big)$.
\end{enumerate}

\begin{claim}
\label{cl1}
For each limit ordinal $\delta<\lambda$,
\[r\forces_{\bbP_\gamma}\mbox{`` }\big(\forall\xi\in w_\delta\big)\big(
\delta\in\mathop{\triangle}\limits_{i<\lambda}\name{A}^\xi_i\big)\
\Rightarrow\ \big(\exists t\in T_\delta\big)\big(\rk_\delta(t)=\gamma\ \&\ 
q^\delta_{*,t}\in \Gamma_{\bbP_\gamma}\big)\mbox{ ''.}\]
\end{claim}

\begin{proof}[Proof of the Claim]
Suppose that $r'\geq r$ and a limit ordinal $\delta<\lambda$ are such that 
\begin{enumerate}
\item[$(*)_{18}$] \quad $r'\forces_{\bbP_\gamma}$`` $\big(\forall\xi\in
  w_\delta\big)\big( \delta\in\mathop{\triangle}\limits_{i<\lambda}
  \name{A}^\xi_i\big)$ ''.  
\end{enumerate}
We are going to show that there is $t\in T_\delta$ such that $\rk_\delta(t)
=\gamma$ and the conditions $q^\delta_{*,t}$ and $r'$ are compatible (and
then the claim will readily follow). To this end let $\langle\vare_\alpha:
\alpha\leq\alpha^*\rangle=w_\delta\cup\{\gamma\}$ be the increasing
enumeration. By induction on $\alpha\leq\alpha^*$ we will choose conditions
$r^*_\alpha,r^{**}_\alpha\in\bbP_{\vare_\alpha}$ and $t=\langle
(t)_{\vare_\alpha}:\alpha<\alpha^*\rangle\in T_\delta$ such that letting  
$t^\alpha_\circ=\langle (t)_{\vare_\beta}:\beta<\alpha\rangle\in T_\delta$
we have  
\begin{enumerate}
\item[$(*)_{19}^\alpha$] $q^\delta_{*,t^\alpha_\circ}\leq r^*_\alpha$ and 
$r'\rest\vare_\alpha\leq r^*_\alpha$,  
\item[$(*)_{20}^\alpha$] $\langle r^*_\beta\conc r'\rest [\vare_\beta,
  \gamma),r^{**}_\beta\conc r'\rest [\vare_\beta,\gamma):\beta<\alpha
  \rangle$ is a partial legal play of $\Game^\lambda_0(\bbP_\gamma,r')$ in
  which Complete uses her winning strategy $\st(\gamma,r')$.    
\end{enumerate}
Suppose that $\alpha\leq\alpha^*$ is a limit ordinal and we have already
defined $t^\alpha_\circ=\langle (t)_{\vare_\beta}: \beta<\alpha\rangle$ and 
$\langle r^*_\beta, r^{**}_\beta: \beta<\alpha\rangle$. Let $\xi=\sup(
\vare_\beta:\beta<\alpha)$. It follows  from $(*)_{20}^\beta$ (for $\beta<
\alpha$) that we may find a condition $s\in \bbP_\xi$ stronger than all
$r^{**}_\beta$ (for $\beta<\alpha$). Let $r^*_\alpha\in\bbP_{\vare_\alpha}$
be such that $r^*_\alpha\rest\xi=s$ and $r^*_\alpha\rest [\xi,\vare_\alpha)
=r'\rest [\xi,\vare_\alpha)$. It follows from $(*)_{19}^\beta$ that 
$q^\delta_{*,t^\alpha_\circ}\rest\xi\leq s=r^*_\alpha\rest\xi$ and $r'\rest
\xi\leq s=r^*_\alpha\rest \xi$. Note also that $(\forall\beta<\alpha)(
r^{**}_\beta\leq s\rest\vare_\beta=r^*_\alpha\rest\vare_\beta)$, so 
$(\forall\beta<\alpha)(r^{**}_\beta\conc r'\rest [\vare_\beta,\gamma)\leq 
r^*_\alpha\conc r'\rest [\vare_\alpha,\gamma))$. Now by induction on $\zeta
\leq \vare_\alpha$ we show that $q^\delta_{*,t^\alpha_\circ}\rest\zeta\leq
r^*_\alpha\rest\zeta$ and $r'\rest\zeta\leq r^*_\alpha\rest\zeta$. For
$\zeta\leq\xi$ we are already done, so assume that $\zeta\in [\xi,
\vare_\alpha)$ and we have shown $q^\delta_{*,t^\alpha_\circ}\rest\zeta \leq
r^*_\alpha\rest\zeta$ and $r'\rest\zeta\leq r^*_\alpha\rest\zeta$. It
follows from $(*)_6+(*)_3$ that $r^*_\alpha\rest\zeta\forces (\forall
i<\delta)(r_i(\zeta)\leq p^\delta_{*,t^\alpha_\circ}(\zeta))$ and therefore
we may use $(*)_{12}$ to conclude that  
\[r^*_\alpha\rest\zeta\forces_{\bbP_\zeta} q^\delta_{*,t^\alpha_\circ}
(\zeta)\leq r_\delta(\zeta)\leq r(\zeta)\leq r'(\zeta)=r^*_\alpha(\zeta).\]  
The limit stages are trivial and we see that $(*)_{19}^\alpha$ and (a part
of) $(*)_{20}^\alpha$ hold. Finally we let $r^{**}_\alpha\in
\bbP_{\vare_\alpha}$ be the condition given to Complete by $\st(\gamma,r')$
as the response to $\langle r^*_\beta\conc r'\rest [\vare_\beta,\gamma),
r^{**}_\beta\conc r'\rest [\vare_\beta,\gamma):\beta<\alpha\rangle\conc
\langle r^*_\alpha\rangle$. 

Now suppose that $\alpha=\beta+1\leq\alpha^*$ and we have already defined
$r^*_\beta, r^{**}_\beta\in\bbP_{\vare_\beta}$ and $t^\beta_\circ\in
T_\delta$. It follows from $(*)_{17}^{\vare_\beta}+(*)_{18}+
(*)_{19}^\beta+(*)_{10}$ that   
\[r^{**}_\beta\forces_{\bbP_{\vare_\beta}}\mbox{`` }r'(\vare_\beta)
\forces_{\name{\bbQ}_{\vare_\beta}}\big(\exists\vare<
\vare^{t^\beta_\circ}_{\delta,\vare_\beta}\big)\big(q^\delta_{*,
t^\beta_\circ\conc\langle\vare\rangle}(\vare_\beta)\in
\Gamma_{\name{\bbQ}_{\vare_\beta}}\big)\mbox{ ''.}\]  
Therefore we may choose $\vare=(t)_{\vare_\beta}<
\vare^{t^\beta_\circ}_{\delta,\vare_\beta}$ (thus defining $t^\alpha_\circ$)
and a condition $s\in\bbP_{\vare_\beta+1}$ such that $s\rest\vare_\beta\geq 
r^{**}_\beta\geq q^\delta_{*, t^\beta_\circ}$ and  
\[s\rest\vare_\beta\forces\mbox{`` }s(\vare_\beta)\geq
r'(\vare_\beta)\ \&\ s(\vare_\beta)\geq q^\delta_{*,t^\alpha_\circ}(
\vare_\beta)\mbox{ ''.}\]
We let $r^*_\alpha\in\bbP_{\vare_\alpha}$ be such that $r^*_\alpha\rest
(\vare_\beta+1)=s$ and $r^*_\alpha\rest (\vare_\beta,\vare_\alpha)=
r'\rest [\vare_\beta,\vare_\alpha)$. Exactly like in the limit case we argue
that $(*)_{19}^\alpha$ and (a part of) $(*)_{20}^\alpha$ hold and then in
the same manner as there we define $r^{**}_\alpha$.  
\medskip

Finally note that $t\in T_\delta$, $\rk_\delta(t)=\gamma$, and the condition 
$r^*_{\alpha^*}$ witnesses that $r'$ and $q^\delta_{*,t}$ are compatible.   
\end{proof}

Now note that 
\[\forces_{\bbP_\gamma}\mbox{`` }\big\{\delta<\lambda:\big(\forall\xi\in
w_\delta\big)\big(\delta\in\mathop{\triangle}\limits_{i<\lambda}
\name{A}^\xi_i\big)\big\}\in\cU^{\bbP_\gamma}\mbox{ '',}\]
and hence by \ref{cl1} we have 
\[r\forces_{\bbP_\gamma}\mbox{`` }\big\{\delta<\lambda:\big(\exists t\in
T_\delta\big)\big(\rk_\delta(t)=\gamma\ \&\ q^\delta_{*,t}\in
\Gamma_{\bbP_\gamma}\big)\big\}\in\cU^{\bbP_\gamma}\mbox{ ''.}\] 
Therefore, by $(*)_7$, 
\[r\forces_{\bbP_\gamma}\mbox{`` }\big\{\delta<\lambda: \big(\exists
\zeta<\zeta_\delta\big)\big(q^\delta_\zeta\in\Gamma_{\bbP_\gamma}\big)
\big\}\in \cU^{\bbP_\gamma}\mbox{ ''}\]
and the proof of the theorem is complete. 
\end{proof}

\begin{remark}
The reason for the weaker ``{\bf b}--bounding'' in the conclusion of
\ref{verB} (and not ``B--bounding'') is that in our description of the
strategy $\st$, we would have to make sure that the conditions played by 
Antigeneric form a tree of conditions. Playing a subgame and keeping the
demands of $(*)_{15}$ are a convenient way to deal with this issue. Note
that (at a stage $\delta$) after playing $\zeta_\delta$ steps of the
subgame, the players may start over and play another $\zeta_\delta$
steps. This small modification can be used to strengthen \ref{verB} to
\ref{verBbis} below.  
\end{remark}

\begin{definition}
\label{doubleres}
Let $\bbQ$ be a strategically $({<}\lambda)$--complete forcing notion.
\begin{enumerate}
\item For a condition $p\in\bbQ$ we define a game $\tbgame(p,\bbQ)$ between 
Generic and Antigeneric as follows. A play of $\tbgame(p,\bbQ)$ lasts
$\lambda$ steps and during a play a sequence      
\begin{enumerate}
\item[$(\boxdot)$]\qquad $\Big\langle \zeta_\alpha,i_\alpha,\langle
p^\alpha_\xi,q^\alpha_\xi:\xi<\zeta_\alpha\cdot (1+i_\alpha)\rangle:\alpha< 
\lambda\Big\rangle$ 
\end{enumerate}
is constructed. Suppose that the players have arrived to a stage $\alpha<
\lambda$ of the game. First, Generic chooses a non-zero ordinal
$\zeta_\alpha<\mu_\alpha$ and then the two players play the first
$\zeta_\alpha$ steps of a subgame in which they alternatively choose
successive terms of a sequence $\langle p^\alpha_\xi,q^\alpha_\xi:\xi<
\zeta_\alpha\rangle$. At a stage $\xi<\zeta_\alpha$ of the subgame, first
Generic picks a condition $p^\alpha_\xi\in\bbQ$ and then Antigeneric answers
with a condition $q^\alpha_\xi$ stronger than $p^\alpha_\xi$. After this
part of the subgame Antigeneric picks a non-zero ordinal $i_\alpha<\lambda$
and the two players continue playing the subgame up to the total length of    
$\zeta_\alpha\cdot (1+i_\alpha)$ alternatively choosing successive terms of a
sequence $\langle p^\alpha_\xi,q^\alpha_\xi:\xi<\zeta_\alpha\cdot
(1+i_\alpha)\rangle$. At a stage $\xi=\zeta_\alpha\cdot i+j$ (where
$j<\zeta_\alpha$, $0<i<1+i_\alpha$) of the subgame, first Generic picks a
$\leq_\bbQ$--upper bound $p^\alpha_\xi\in\bbQ$ to $\{q^\alpha_\xi:\xi=
\zeta_\alpha\cdot i'+j\ \&\ i'<i\}$, and then Antigeneric answers with a
condition $q^\alpha_\xi$ stronger than $p^\alpha_\xi$.  At the end, Generic
wins the play  $(\boxdot)$ of $\tbgame(p,\bbQ)$ if and only if both players
had always legal moves and   
\begin{enumerate}
\item[$(\circledast)^{\rm rc}_{\bf 2b}$] there is a condition $p^*\in\bbQ$ 
stronger than $p$ and such that  
\[p^*\forces_{\bbQ}\mbox{`` }\big\{\alpha<\lambda:\big(\exists j<
\zeta_\alpha\big)\big(\forall i<1+i_\alpha\big)\big(q^\alpha_{\zeta_\alpha
  \cdot i+j}\in\Gamma_{\bbQ}\big)\big\} \in \cU^\bbQ \mbox{ ''}.\] 
\end{enumerate}
\item The game $\tagame(p,\bbQ)$ is defined similarly except that the
  winning criterion $(\circledast)^{\rm rc}_{\bf 2b}$ is 
  changed so that ``$\in \cU^\bbQ$'' is replaced by ``$=\lambda$''.
\item We say that a forcing notion $\bbQ$ is {\em reasonably double {\bf
    a}--bounding over $\bar{\mu}$\/} if    
\begin{enumerate}
\item[(a)] $\bbQ$ is strategically $({<}\lambda)$--complete,  and 
\item[(b)] for any $p\in\bbQ$, Generic has a winning strategy in the game
$\tagame(p,\bbQ)$. 
\end{enumerate}
In an analogous manner we define when the forcing notion $\bbQ$ is {\em
reasonably double {\bf b}--bounding over $\cU,\bar{\mu}$}.
\end{enumerate}
\end{definition}

\begin{observation}
\label{doubleobs}
Let $\cU,\bar{\mu}$ be as in \ref{incon}, $(X,x)\in\{(A,{\bf a}),(B,{\bf
b})\}$.  Then the following implications hold for a forcing notion $\bbQ$:  
\begin{center}
\begin{tabular}{c}
``$\bbQ$ is reasonably $X$--bounding over $\cU,\bar{\mu}$''\\
$\Downarrow$\\
``$\bbQ$ is reasonably double $x$--bounding over $\cU,\bar{\mu}$''\\
$\Downarrow$\\
$\bbQ$ is reasonably $x$--bounding over $\cU,\bar{\mu}$'' .
\end{tabular}
\end{center}
\end{observation}

\begin{corollary}
\label{verBbis}
Assume that $\lambda,\bar{\mu}$ are as in \ref{incon} and $\bar{\bbQ}= 
\langle\bbP_\xi,\name{\bbQ}_\xi:\xi<\gamma\rangle$ is a $\lambda$--support
iteration such that for every $\xi<\gamma$,  
\[\forces_{\bbP_\xi}\mbox{`` $\name{\bbQ}_\xi$ is reasonably B--bounding
  over $\bar{\mu}$ ''.}\] 
Then $\bbP_{\gamma}=\lim(\bar{\bbQ})$ is reasonably double ${\bf
b}$--bounding over $\bar{\mu}$.
\end{corollary}

\begin{theorem}
\label{verA}
Assume that $\lambda,\bar{\mu}$ are as in \ref{incon} and $\bar{\bbQ}= 
\langle\bbP_\xi,\name{\bbQ}_\xi:\xi<\gamma\rangle$ is a $\lambda$--support
iteration such that for every $\xi<\gamma$,  
\[\forces_{\bbP_\xi}\mbox{`` $\name{\bbQ}_\xi$ is reasonably A--bounding
  over $\bar{\mu}$ ''.}\] 
Then $\bbP_{\gamma}=\lim(\bar{\bbQ})$ is reasonably double ${\bf
a}$--bounding over $\bar{\mu}$ (and thus also reasonably ${\bf
a}$--bounding over $\bar{\mu}$.
\end{theorem}

\begin{proof}
The proof of Theorem \ref{verB} (changed so that it works for \ref{verBbis}) 
can be easily modified to fit the current  purpose, just replace each
occurrence of $\Game^{\rm rcB}_{\cU,\bar{\mu}}, \Game^{{\rm rc}{\bf
2b}}_{\cU,\bar{\mu}}$ by $\Game^{\rm rcA}_{\bar{\mu}}, \Game^{{\rm
rc}{\bf 2a}}_{\bar{\mu}}$, respectively, and in the end think 
that $\forces_{\bbP_{\xi+1}}\name{A}^\xi_i=\lambda$ (for $\xi\in\Dom(r)$,
$i<\lambda$). Then one uses the proof of \ref{cl1} to argue that  
\[r\forces_{\bbP_\gamma}\mbox{`` }\big(\forall\delta<\lambda\big)
\big(\exists t\in T_\delta\big)\big(\rk_\delta(t)=\gamma\ \&\
q^\delta_{*,t}\in \Gamma_{\bbP_\gamma}\big)\mbox{ ''.}\] 
\end{proof}

\begin{problem}
\label{prob}
\begin{enumerate}
\item Do we have a result parallel to \ref{verB} and/or \ref{verA} for reasonably
C--bounding forcings? 
\item Can the implications in \ref{doubleobs} be reversed in the sense that
  we allow passing to an equivalent forcing notion? 
\end{enumerate}
\end{problem}

\section{Consequences of reasonable ABC} 

Let us note that Theorem \ref{verA} improves \cite[Theorem
A.2.4]{RoSh:777}. Before we explain why, we should recall the following
definition.  

\begin{definition}
[{\cite[Def. A.2.1]{RoSh:777}}]
\label{da2}
Let $\bbP$ be a forcing notion.
\begin{enumerate}
\item A complete $\lambda$--tree of height $\alpha<\lambda$ is a set of
  sequences $s\subseteq {}^{\leq\alpha}\lambda$ such that 
\begin{itemize}
\item $s$ has the $\vtl$--smallest element denoted $\mrot(s)$, 
\item $s$ is closed under initial segments longer than $\lh(\mrot(s))$,  and 
\item the union of any $\vtl$--increasing sequence of members of $s$ is in
  $s$, and 
\item $\big(\forall\eta\in s\big)\big(\lh(\eta)\leq\alpha\big)$ and
  $\big(\forall\eta\in s\big)\big(\exists\nu\in s\big)\big(\eta\vtl\nu\ \&\
  \lh(\nu)=\alpha\big)$. 
\end{itemize}
\item For a condition $p\in\bbP$ and an ordinal $i_0<\lambda$ we define a
game $\Gsg(i_0,p,\bbP)$ of two players, {\em the Generic player} and {\em
the Antigeneric player}. A play lasts $\lambda$ moves indexed by ordinals
from the interval $[i_0,\lambda)$, and during it the players construct
a sequence $\langle (s_i,\bar{p}^i,\bar{q}^i):i_0\leq i<\lambda\rangle$ as
follows. At stage $i$ of the play (where $i_0\leq i<\lambda$), first 
Generic chooses $s_i\subseteq {}^{\leq i+1}\lambda$ and a system
$\bar{p}^i=\langle p^i_\eta:\eta\in s_i\cap {}^{ i+1}\lambda
\rangle$ such that  
\begin{enumerate}
\item[$(\alpha)$] $s_i$ is a complete $\lambda$--tree of height $i+1$ and 
$\lh(\mrot(s_i))=i_0$, 
\item[$(\beta)$]  for all $j$ such that $i_0\leq j<i$ we have $s_j=s_i\cap
{}^{\leq j+1}\lambda$,
\item[$(\gamma)$] $p^i_\eta\in\bbP$ for all $\eta\in s_i\cap
{}^{ i+1}\lambda$, and 
\item[$(\delta)$] if $i_0\leq j<i$, $\nu\in s_i\cap {}^{ j+1}
\lambda$ and $\nu\vtl\eta\in s_i\cap {}^{ i+1}\lambda$, then
$q^j_\nu\leq p^i_\eta$ and $p\leq p^i_\eta$, 
\item[$(\vare)$] $|s_i\cap {}^{ i+1}\lambda|<\mu_i$.
\end{enumerate}
Then Antigeneric answers choosing a system $\bar{q}^i=\langle q^i_\eta:
\eta\in s_i\cap {}^{ i+1}\lambda\rangle$ of conditions in $\bbP$
such that $p^i_\eta\leq q^i_\eta$ for each $\eta\in s_i\cap {}^{
i+1}\lambda$.  

The Generic player wins a play if she always has legal moves (so the play
really lasts $\lambda$ steps) and there are a condition $q\geq p$ and a
$\bbP$--name $\name{\rho}$ such that 
\begin{enumerate}
\item[$(\circledast)$]\qquad $q\forces_{\bbP}\mbox{`` }\name{\rho}\in
\bairel\ \&\ \big(\forall i\in [i_0,\lambda\big))\big(\name{\rho}\rest (i+1)
\in s_i\ \&\ q^i_{\name{\rho}\rest (i+1)}\in\Gamma_{\bbP}\big)\mbox{ ''.}$
\end{enumerate}
\item We say that $\bbP$ has the {\em strong $\bar{\mu}$--Sacks
property\/} whenever 
\begin{enumerate}
\item[(a)] $\bbP$ is strategically $(<\lambda)$--complete, and 
\item[(b)] the Generic player has a winning strategy in the game
$\Gsg(i_0,p,\bbP)$ for any $i_0<\lambda$ and $p\in\bbP$.
\end{enumerate}
\end{enumerate}
\end{definition}

The following proposition explains why \ref{verA} is stronger than
\cite[Theorem A.2.4]{RoSh:777}.

\begin{proposition}
\label{thesame}
Assume that $\lambda,\bar{\mu}$ are as in Context \ref{incon} and that
additionally $(\forall i<j<\lambda)(\mu_i\leq \mu_j)$. Let $\bbQ$ be a
forcing notion. Then 
\begin{center}
$\bbQ$ is reasonably A--bounding over $\bar{\mu}$
\end{center}
if and only if 
\begin{center}
$\bbQ$ has the strong $\bar{\mu}$--Sacks property.
\end{center}
\end{proposition}

\begin{proof}
Suppose that $\bbQ$ is reasonably A--bounding over $\bar{\mu}$. Since the
sequence $\bar{\mu}$ is non-decreasing, it is enough to show that Generic
has a winning strategy in $\Gsg(0,p,\bbQ)$ for each $p\in\bbQ$ (as then
almost the same strategy will be good in $\Gsg(i,p,\bbQ)$ for any
$i<\lambda$). 

Let $p\in\bbQ$. We are going to define a strategy $\st$ for Generic in the
game $\Gsg(0,p,\bbQ)$. To describe it, let us fix a winning strategy $\st_0$
of Complete in $\Game^\lambda_0(\bbQ,p)$ and a winning strategy $\st_1$ of
Generic in $\Agame(p,\bbQ)$. Now, at a stage $\delta<\lambda$ of the play
the strategy $\st$ will tell Generic to write aside  
\begin{enumerate}
\item[$(\boxtimes)_\delta$]\quad $I_\delta$ and $\langle r^{0,\delta}_t,
r^{1,\delta}_t:t\in I_\delta\rangle$ and $\langle r^\delta_\eta: \eta\in
s_\delta\cap {}^{\delta+1}\lambda\rangle$ 
\end{enumerate}
so that if $\langle (s_\delta,\bar{p}^\delta,\bar{q}^\delta):\delta<\lambda
\rangle$ is a play of $\Gsg(0,p,\bbQ)$  in which Generic follows $\st$, then
the following conditions $(\odot)_1$--$(\odot)_4$ are satisfied (for each
$\delta<\lambda$).     
\begin{enumerate}
\item[$(\odot)_1$] $\big\langle I_\alpha,\langle r^{0,\alpha}_t,
r^{1,\alpha}_t: t\in I_\alpha\rangle:\alpha\leq\delta\big\rangle$ is a
partial legal play of $\Agame(p,\bbQ)$ in which Generic uses $\st_1$. 
\item[$(\odot)_2$] For each $\eta\in s_\delta\cap {}^{\delta+1}\lambda$ the 
sequence $\langle q^\alpha_{\eta\rest (\alpha+1)},r^\alpha_{\eta\rest
(\alpha+1)}:\alpha\leq\delta\rangle$ is a partial legal play of
$\Game^\lambda_0 (\bbQ,p)$ in which Complete uses $\st_0$. 
\item[$(\odot)_3$] If $t\in I_\delta$, $\alpha<\delta$, $\nu\in s_\alpha\cap 
  {}^{\alpha+1}\lambda$, then either $r^\alpha_\nu,r^{1,\delta}_t$ are
  incompatible or $r^\alpha_\nu\leq r^{1,\delta}_t$.
\item[$(\odot)_4$] $\langle p^\delta_\nu:\nu\in s_\delta\cap{}^{\delta+1}
\lambda \rangle$ is an antichain in $\bbQ$. 
\end{enumerate}
So suppose that the two players arrived to a stage $\delta<\lambda$ of the
game $\Gsg(0,p,\bbQ)$ and the objects listed in $(\boxtimes)_\alpha$ (for
$\alpha<\delta$) as well as $\langle (s_\alpha,\bar{p}^\alpha,
\bar{q}^\alpha):\alpha<\delta\rangle$ have been constructed. First Generic
uses $\st_1$ to pick the answer $\big(I_\delta,\langle r^{0,\delta}_t:t\in
I_\delta\rangle\big)$ to $\big\langle I_\alpha,\langle
r^{0,\alpha}_t,r^{1,\alpha}_t:t\in I_\alpha\rangle:\alpha<\delta\big\rangle$  
in $\Agame(p,\bbQ)$. Then she uses the strategic completeness of $\bbQ$ and
\ref{obsA.4} to choose a system $\langle r^*_t:t\in I_\delta\rangle$ of
conditions in $\bbQ$ such that  
\begin{enumerate}
\item[$(\odot)_5$] if $t\in I_\delta$, then $r^{0,\delta}_t\leq r^*_t$ and
for every $\alpha<\delta$ and $\nu\in s_\alpha\cap {}^{\alpha+1}\lambda$,
either $r^\alpha_\nu,r^*_t$ are incompatible or $r^\alpha_\nu\leq r^*_t$,
and also either $p,r^*_t$ are incompatible or $p\leq r^*_t$,
\item[$(\odot)_6$] if $t_0,t_1\in I_\delta$, $t_0\neq t_1$, then the
  conditions $r^*_{t_0},r^*_{t_1}$ are incompatible in $\bbQ$. 
\end{enumerate}
Now she lets $s^*=\{\eta\in {}^\delta\lambda:(\forall\alpha<\delta)(\eta
\rest (\alpha+1)\in s_\alpha)\}$ and
\[s^-=\{\eta\in s^*: (\exists t\in I_\delta)(\forall\alpha<\delta)(
r^\alpha_{\eta\rest (\alpha+1)}\leq r^*_t\ \&\ p\leq r^*_t)\},\]
and for each $\eta\in s^-$ she fixes an enumeration $\langle t^\eta_\xi:
\xi<\xi_\eta\rangle$ of the set 
\[\big\{t\in I_\delta:\big(\forall\alpha<\delta\big)\big(r^\alpha_{\eta
\rest (\alpha+1)}\leq r^*_t\ \&\ p\leq r^*_t\big)\big\}.\] 
Now Generic defines
\[s_\delta^+=\big\{\nu\in {}^{\delta+1}\lambda:\big(\nu\rest\delta\in s^*
\setminus s^-\ \&\ \nu(\delta)=0\big)\mbox{ or }\big(\nu\rest\delta\in s^-\
\&\ \nu(\delta)<\xi_{\nu\rest\delta}\big)\big\}\]
and she lets $s_\delta$ be a $\lambda$--tree of height $\delta+1$ such that
$s_\delta\cap {}^{\delta+1}\lambda=s^+_\delta$. For $\nu\in s^+_\delta$ she
also chooses $p^\delta_\nu$ so that 
\begin{itemize}
\item if $\nu\rest\delta\notin s^-$, then $p^\delta_\nu\in\bbQ$ is an 
upper bound to $\{r^\alpha_{\nu\rest (\alpha+1)}:\alpha<\delta\}\cup\{p\}$
(remember $(\odot)_2$), 
\item if $\nu\rest\delta\in s^-$, then $p^\delta_\nu= r^*_{t^{\nu\rest
\delta}_{\nu(\delta)}}$.
\end{itemize}
And now, in the play of $\Gsg(0,p,\bbQ)$, Generic puts 
\[s_\delta\quad\mbox{ and }\quad \langle p^\delta_\nu:\nu\in s^+_\delta
\rangle\]
and Antigeneric answers with $\langle q^\delta_\nu:\nu\in s^+_\delta\rangle$
(so that $q^\delta_\nu\geq p^\delta_\nu$). Conditions $r^\delta_\nu$ (for
$\nu\in s^+_\delta$) are determined using $\st_0$ (so that the demand in
$(\odot)_2$ is satisfied). Finally, Generic defines also $r^{1,\delta}_t$
for $t\in I_\delta$ so that  
\begin{itemize}
\item if $t=t^\eta_\xi$ for some $\eta\in s^-$ and $\xi<\xi_\eta$, then
  $r^{1,\delta}_t=r^\delta_{\eta\conc\langle\xi\rangle}$,  
\item otherwise $r^{1,\delta}_t=r^*_t$.
\end{itemize}
This completes the description of what Generic plays and what she writes
aside --- it should be clear that the requirements of
$(\odot)_1$--$(\odot)_4$ are satisfied. Now, why is $\st$ a winning
strategy? So suppose that $\langle (s_\delta,\bar{p}^\delta,\bar{q}^\delta):
\delta<\lambda\rangle$ is a play of $\Gsg(0,p,\bbQ)$ in which Generic
follows $\st$, and $I_\delta$, $\langle r^{0,\delta}_t,r^{1,\delta}_t:t\in
I_\delta\rangle$ and $\langle r^\delta_\eta: \eta\in s_\delta\cap
{}^{\delta+1}\lambda\rangle$ (for $\delta<\lambda$) are the objects
constructed by Generic aside, so they satisfy $(\odot)_1$--$(\odot)_4$. 
It follows from $(\odot)_1$ and the choice of $\st_1$ that there is a
condition $p^*\geq p$ such that 
\begin{enumerate}
\item[$(\odot)_7$] for every $\delta<\lambda$ the set $\big\{r^{1,\delta}_t:
  t\in I_\delta\big\}$ is pre-dense above $p^*$. 
\end{enumerate}
We claim that then also 
\begin{enumerate}
\item[$(\odot)_8$] for every $\delta<\lambda$ the set $\big\{r^\delta_\eta:
  \eta\in s_\delta\cap {}^{\delta+1}\lambda\big\}$ is pre-dense above $p^*$ 
\end{enumerate}
(and this clearly implies that Generic won the play, remember
$(\odot)_4$). Assume towards contradiction that $(\odot)_8$ fails and let
$\delta<\lambda$ be the smallest ordinal for which we may find a condition
$q\geq p^*$ such that $q$ is incompatible with every $r^\delta_\eta$ for
$\eta\in s_\delta\cap {}^{\delta+1}\lambda$. It follows from $(\odot)_7$
that we may pick $t\in I_\delta$ such that the conditions $r^{1,\delta}_t,q$
are compatible. By the previous sentence and by the definition of 
$r^{1,\delta}_t$ we get that $t\neq t^\eta_\xi$ for all $\xi<\xi_\eta$,
$\eta\in s^-$ and thus $r^{1,\delta}_t=r^*_t$. Look at the condition $r^*_t$
(satisfying $(\odot)_5+(\odot)_6$) --- it must be stronger than $p$ and by
the minimality of $\delta$ we have that $\big(\forall\alpha<\delta\big)
\big(\exists\nu\in s_\alpha\cap {}^{\alpha+1}\lambda\big)\big(r^\alpha_\nu
\leq r^*_t\big)$. It follows from $(\odot)_4$ from stages $\alpha<\delta$
that there is $\eta\in s^*$ such that $\big(\forall\alpha<\delta\big)
\big(r^\alpha_{\eta\rest(\alpha+1)}\leq r^*_t\big)$. Then $t\in s^-$ and
hence $t=t^\eta_\xi$ for some $\xi<\xi_\eta$, contradicting what we already
got. 
\medskip

The converse implication should be clear.
\end{proof}

The following easy proposition explains why the names of the properties
defined in \ref{p.1A} include the adjective ``bounding''. 

\begin{proposition}
\label{bound}
Let $\lambda,\cU$ and $\bar{\mu}$ be as in \ref{incon}. Assume that $\bbQ$
is a forcing notion, $p\in \bbQ$ and $\name{\tau}$ is a $\bbQ$--name for an
element of $\bairel$. 
\begin{enumerate}
\item If $\bbQ$ is reasonably {\bf a}--bounding over $\bar{\mu}$, then there
  are a condition $q\geq p$ and a sequence $\bar{a}=\langle a_\alpha:\alpha<
  \lambda\rangle$ such that 
\begin{itemize}
\item $a_\alpha\subseteq\lambda$, $|a_\alpha|<\mu_\alpha$ for all $\alpha<
  \lambda$, 
\item $q\forces_{\bbQ}$`` $(\forall\alpha<\lambda)(\name{\tau}(\alpha)\in
  a_\alpha)$ ''.
\end{itemize}
\item If $\bbQ$ is reasonably {\bf b}--bounding over $\cU,\bar{\mu}$, then
  there are a condition $q\geq p$ and a sequence $\bar{a}=\langle
  a_\alpha:\alpha<\lambda\rangle$ such that 
\begin{itemize}
\item $a_\alpha\subseteq\lambda$, $|a_\alpha|<\mu_\alpha$ for all $\alpha<
  \lambda$, 
\item $q\forces_{\bbQ}$`` $\{\alpha<\lambda:\name{\tau}(\alpha)\in
  a_\alpha\}\in\cU^{\bbQ}$ ''. 
\end{itemize}
\item If $\bbQ$ is reasonably {\bf c}--bounding over $\cU,\bar{\mu}$, then
  there are a condition $q\geq p$ and a sequence $\bar{a}=\langle
  a_\alpha:\alpha<\lambda\rangle$ such that 
\begin{itemize}
\item $a_\alpha\subseteq\lambda$, $|a_\alpha|<\mu_\alpha$ for all $\alpha<
  \lambda$, 
\item $q\forces_{\bbQ}$`` $\{\alpha<\lambda:\name{\tau}(\alpha)\in
  a_\alpha\}\in\big(\cU^{\bbQ}\big)^+$ ''. 
\end{itemize}
\end{enumerate}
\end{proposition}

\section{Forcing notions and models}
In this section, in  addition to the assumptions stated in \ref{incon}
we will also assume that 

\begin{context}
\label{extraassum}
\begin{enumerate}
\item[(d)] $S\subseteq \lambda$ is stationary and co-stationary, $S\in\cU$,   
\item[(e)] $\cV$ is a normal filter on $\lambda$, $\lambda\setminus
  S\in\cV$.     
\end{enumerate}
\end{context}

\begin{definition}
\label{ext}
\begin{enumerate}
\item Let $\alpha<\beta<\lambda$. {\em An $(\alpha,\beta)$--extending
function\/} is a mapping $c:\cP(\alpha)\longrightarrow\cP(\beta)\setminus 
\cP(\alpha)$ such that $c(u)\cap\alpha=u$ for all $u\in\cP(\alpha)$. 
\item Let $C$ be an unbounded subset of $\lambda$. {\em A $C$--extending
sequence\/} is a sequence $\gc=\langle c_\alpha:\alpha\in C\rangle$ such
that each $c_\alpha$ is an $(\alpha,\min(C\setminus(\alpha+1)))$--extending
function.  
\item Let $C\subseteq\lambda$, $\|C\|=\lambda$, $\beta\in C$, $w\subseteq
\beta$ and let $\gc=\langle c_\alpha:\alpha\in C\rangle$ be a $C$--extending
sequence. We define $\pos^+(w,\gc,\beta)$ as the family of all subsets $u$
of $\beta$ such that  
\begin{enumerate}
\item[(i)]  if $\alpha_0=\min\big(\{\alpha\in C:(\forall\xi\in w)(\xi<
\alpha)\}\big)$, then $u\cap\alpha_0=w$ (so if $\alpha_0=\beta$, then
$u=w$), and  
\item[(ii)] if $\alpha_0,\alpha_1\in C$, $w\subseteq\alpha_0<\alpha_1=
\min(C\setminus(\alpha_0+1))\leq\beta$, then either $c_{\alpha_0}(u\cap
\alpha_0)=u\cap\alpha_1$ or $u\cap\alpha_0=u\cap\alpha_1$, 
\item[(iii)] if $\sup(w)<\alpha_0=\sup(C\cap\alpha_0)\notin C$,
$\alpha_1=\min\big(C\setminus(\alpha_0+1)\big)\leq\beta$, then
$u\cap\alpha_1=u\cap\alpha_0$.
\end{enumerate}
The family $\pos(w,\gc,\beta)$ consists of all elements $u$ of
$\pos^+(w,\gc,\beta)$ which satisfy also the following condition:
\begin{enumerate}
\item[(iv)] if $\alpha_0=\min\big(\{\alpha\in C:w\subseteq\alpha\}\big)\leq
\beta$, $\alpha_1=\min\big(C\setminus(\alpha_0+1)\big)\leq\beta$, then $u
\cap\alpha_1=c_{\alpha_0}(w)$.
\end{enumerate}
\item A $C$--extending sequence $\gc=\langle c_\alpha:\alpha\in C\rangle$ is
{\em $S$--closed\/} provided that 
\begin{enumerate}
\item[(i)]   $C$ is a club of $\lambda$, and 
\item[(ii)]  if $\alpha\in C$ and $u\subseteq\alpha$, then $\alpha\in
  c_\alpha(u)$, and 
\item[(iii)] if $\xi\in S\setminus C$, $\alpha\in C\cap\xi$, $u\subseteq
\alpha$ and $\xi=\sup\big(c_\alpha(u)\cap\xi\big)$, then $\xi\in
c_\alpha(u)$.  
\end{enumerate}
\item A set $w\subseteq\lambda$ is $S$--closed if $\xi=\sup\big(w\cap\xi
  \big)\in S$ implies $\xi\in w$.
\item Let $\gc=\langle c_\alpha:\alpha\in C\rangle$ be an $S$--closed
$C$--extending sequence, $\beta\in C$, $w\subseteq\beta$ and $\alpha=\min
\big(C\setminus\sup(w)\big)$. Assume also that $w\cup\{\alpha\}$ is
$S$--closed. Then we let 
\[\begin{array}{ll}
\pos^+_S(w,\gc,\beta)=&\big\{u\in\pos^+(w,\gc,\beta): u\cup\{\beta\}\mbox{
  is $S$--closed }\big\}\\ 
\pos_S(w,\gc,\beta)=&\big\{u\in\pos(w,\gc,\beta): u\cup\{\beta\}\mbox{ is
  $S$--closed }\big\}
  \end{array}\]
\end{enumerate}
\end{definition}

\begin{observation}
\label{easyobs}
\begin{enumerate}
\item Assume that $\gc$ is a $C$--extending sequence, $\alpha,\beta\in C$,
$\alpha<\beta$ and $w\subseteq\alpha$.
\begin{enumerate}
\item If $u\in\pos(w,\gc,\alpha)$ and $v\in\pos(u,\gc,\beta)$, then $v\in
\pos(w,\gc,\beta)$. 
\item If $v\in\pos(w,\gc,\beta)$, then $v\cap\alpha\in\pos(w,\gc,\alpha)$
and $v\in\pos^+(v\cap\alpha,\gc,\beta)$.  
\item Similarly for $\pos^+$.
\end{enumerate}
\item Assume that $\gc$ is an $S$--closed  $C$--extending sequence,
$\alpha,\beta\in C$, $\alpha<\beta$, $w\subseteq\alpha$ and $w\cup\big\{
\min\big(C\setminus\sup(w)\big)\big\}$ is $S$--closed. 
\begin{enumerate}
\item If $u\in\pos_S(w,\gc,\alpha)$ and $v\in\pos_S(u,\gc,\beta)$, then
  $v\in \pos_S(w,\gc,\beta)$. 
\item If $v\in\pos_S(w,\gc,\beta)$, then $v\cap\alpha\in\pos_S(w,\gc,\alpha)$ 
and $v\in\pos^+_S(v\cap\alpha,\gc,\beta)$.  
\item Similarly for $\pos^+_S$.
\item $\emptyset\neq \pos_S(w,\gc,\beta)\subseteq\pos^+_S(w,\gc,\beta)$.
\end{enumerate}
\end{enumerate}
\end{observation}

\begin{definition}
\label{clfor}
We define a forcing notion $\clfo$ as follows.\\
{\bf A condition in $\clfo$} is a triple $p=(w^p,C^p,\gc^p)$ such that 
\begin{enumerate}
\item[(i)]    $C^p\subseteq\lambda$ is a club of $\lambda$ and $w^p\subseteq 
  \min(C^p)$ is such that the set $w^p\cup\{\min(C^p)\}$ is $S$--closed,  
\item[(ii)]   $\gc^p=\langle c^p_\alpha:\alpha\in C^p\rangle$ is an
  $S$--closed $C^p$--extending sequence.
\end{enumerate}
{\bf The order $\leq_{\clfo}=\leq$ of $\clfo$} is given by \\
$p\leq_{\clfo} q$\qquad if and only if \\
\begin{enumerate}
\item[(a)] $C^q\subseteq C^p$ and $w^q\in\pos^+_S(w^p,\gc^p,\min(C^q))$ and  
\item[(b)] if $\alpha_0<\alpha_1$ are two successive members of $C^q$, $u\in
\pos^+_S(w^q,\gc^q,\alpha_0)$, then $c^q_{\alpha_0}(u)\in\pos_S(u,\gc^p,
\alpha_1)$. 
\end{enumerate}
For $p\in\clfo$, $\alpha\in C^p$ and $u\in\pos^+_S(w^p,\gc^p,\alpha)$ we let  
$p\rest_\alpha u\stackrel{\rm def}{=}(u,C^p\setminus\alpha,\gc^p\rest(C^p
\setminus\alpha))$.
\end{definition}

\begin{proposition}
\label{clbas}
\begin{enumerate}
\item $\clfo$ is a $({<}\lambda)$--complete forcing notion of cardinality 
  $2^\lambda$. 
\item If $p\in\clfo$ and $\alpha\in C^p$, then
\begin{itemize}
\item for each $u\in\pos^+_S(w^p,\gc^p,\alpha)$, $p\rest_\alpha u\in\clfo$ is
  a condition stronger than $p$, and    
\item the family $\{p\rest_\alpha u: u\in\pos^+_S(w^p,\gc^p,\alpha)\}$ is
  pre-dense above $p$. 
\end{itemize}
\item Let $p\in\clfo$ and $\alpha<\beta$ be two successive members of
$C^p$. Suppose that for each $u\in\pos^+_S(w^p,\gc^p,\alpha)$ we are given a
condition $q_u\in\clfo$ such that $p\rest_\beta c^p_\alpha(u)\leq q_u$. 
Then there is a condition $q\in\clfo$ such that letting $\alpha'=\min(
C^q\setminus\beta)$ we have
\begin{enumerate}
\item[(a)] $p\leq q$, $w^q=w^p$, $C^q\cap\beta=C^p\cap\beta$ and
$c^q_\delta=c^p_\delta$ for $\delta\in C^q\cap\alpha$, and  
\item[(b)] $\bigcup\big\{w^{q_u}:u\in\pos^+_S(w^p,\gc^p,\alpha)\big\}
\subseteq\alpha'$, and 
\item[(c)] $q_u\leq q\rest_{\alpha'} c^q_\alpha(u)$ for every $u\in\pos^+_S(
w^p,\gc^p,\alpha)$. 
\end{enumerate}
\item Assume that $p\in\clfo$, $\alpha\in C^p$ and $\name{\tau}$ is a
$\clfo$--name such that $p\forces$``$\name{\tau}\in\bV$''. Then there is a
  condition $q\in\clfo$ stronger than $p$ and such that   
\begin{enumerate}
\item[(a)] $w^q=w^p$, $\alpha\in C^q$ and $C^q\cap\alpha=C^p\cap\alpha$,
  and 
\item[(b)] if $u\in\pos^+_S(w^q,\gc^q,\alpha)$ and $\gamma=\min(C^q\setminus 
(\alpha+1))$, then the condition $q\rest_\gamma c^q(u)$ forces a value to
$\name{\tau}$. 
\end{enumerate}
\end{enumerate}
\end{proposition}

\begin{proof}  
(1)\quad It should be clear that $\clfo$ is a forcing notion of size
$2^\lambda$. To show that it is $({<}\lambda)$--complete suppose that
$\gamma<\lambda$ is a limit ordinal and $\bar{p}=\langle p_\xi:\xi<
\gamma\rangle\subseteq\clfo$ is $\leq_{\clfo}$--increasing. We put 
$w^q=\bigcup\limits_{\xi<\gamma}w^{p_\xi}$, $C^q=\bigcap\limits_{\xi<
\gamma} C^{p_\xi}$ and for $\delta\in C^q$ we define $c^q_\delta:
\cP(\delta)\longrightarrow\cP\big(\min(C^q\setminus(\delta+1))\big)$ so that  
\begin{itemize}
\item if $u\in\bigcap\limits_{\xi<\gamma}\pos^+_S(w^{p_\xi},\gc^{p_\xi},
\delta)$, then $c^q_\delta(u)=\bigcup\limits_{\xi<\gamma}c^{p_\xi}_\delta(
u)$, 
\item if $u\subseteq\delta$ but it is not in $\bigcap\limits_{\xi<\gamma}
\pos^+_S(w^{p_\xi},\gc^{p_\xi},\delta)$, then $c^q_\delta(u)=u\cup
\{\delta\}$.  
\end{itemize}
Finally we put $\gc^q=\langle c^q_\delta:\delta\in
C^q\rangle$ and $q=(w^q,C^q,\gc^q)$. One easily checks that $q\in\clfo$ is a 
condition stronger than all $p_\xi$'s. 
\medskip

\noindent (2)\quad Straightforward (remember \ref{easyobs}(2)).
\medskip

\noindent (3)\quad We let $w^q=w^p$ and $C^q=(C^p\cap\beta)\cup\bigcap\big\{
C^{q_u}:u\in\pos^+_S(w^p,\gc^p,\alpha)\big\}$ (plainly, $C^q$ is a club of
$\lambda$). Let $\alpha'=\min(C^q\setminus(\alpha+1))=\min(C^q\setminus
\beta)$. For $\delta\in C^q\cap\alpha=C^p\cap\alpha$ put $c^q_\delta=
c^p_\delta$. Next, choose an $(\alpha,\alpha')$--extending function
$c^q_\alpha:\cP(\alpha)\longrightarrow\cP(\alpha')$ such that $(\forall u\in
\pos^+_S(w^p,\gc,\alpha))(c^q_\alpha(u)\in\pos_S(w^{q_u},\gc^{q_u},\alpha'))$
and $(c^q_\alpha(u)\setminus\alpha)\cup\{\alpha'\}$ is $S$--closed for each
$u\subseteq\alpha$. (Remember \ref{easyobs}(2d); note that, by the
definition of $C^q$, $w^{q_u}\subseteq\alpha'$ for each $u\in\pos^+_S(w^p,
\gc^p,\alpha)$.) Finally, if $\delta_0<\delta_1$ are two successive members
of $C^q\setminus\alpha'$, then choose a $(\delta_0,\delta_1)$--extending
function  $c^q_{\delta_0}:\cP(\delta_0)\longrightarrow\cP(\delta_1)$ so that  
\begin{enumerate}
\item[(i)]  if $v\subseteq\delta_0$, $u=v\cap\alpha\in\pos^+_S(w^p,\gc^p,
\alpha)$ and $v\in\pos^+_S(w^{q_u},\gc^{q_u},\delta_0)$, then
$c^q_{\delta_0}(v)\in\pos_S(v,\gc^{q_u},\delta_1)$;
\item[(ii)] if $v\subseteq\delta_0$ but we are not in a case covered by
  (i), then $c^q_{\delta_0}(v)\in\pos_S(v,\gc^p,\delta_1)$.
\end{enumerate}
Let $\gc^q=\langle c^q_\delta:\delta\in C^q\rangle$ and
$q=(w^q,C^q,\gc^q)$. It should be clear that $q\in\clfo$ is a condition as
required.   
\medskip

\noindent (4)\quad Easily follows from (3).
\end{proof}

\begin{definition}
\label{natural}
Suppose that $\gamma<\lambda$ is a limit ordinal and $\bar{p}=\langle
p_\xi:\xi<\lambda\rangle\subseteq\clfo$ is $\leq_{\clfo}$--increasing. The
condition $q$ constructed as in the proof of \ref{clbas}(1) for $\bar{p}$
will be called {\em the natural limit of $\bar{p}$}. 
\end{definition}

\begin{proposition}
\label{cllimit}
\begin{enumerate}
\item Suppose $\bar{p}=\langle p_\xi:\xi<\lambda\rangle$ is a
$\leq_{\clfo}$--increasing sequence of conditions from $\clfo$ such that 
\begin{enumerate}
\item[(a)] $w^{p_\xi}=w^{p_0}$ for all $\xi<\lambda$, and 
\item[(b)] if $\gamma<\lambda$ is limit, then $p_\gamma$ is the natural
  limit of $\bar{p}\rest\gamma$, and 
\item[(c)] for each $\xi<\lambda$, if $\delta\in C^{p_\xi}$, $\otp(C^{p_\xi}
\cap\delta)=\xi$, then $C^{p_{\xi+1}}\cap (\delta+1)=C^{p_\xi}\cap
(\delta+1)$ and for every $\alpha\in C^{p_{\xi+1}}\cap\delta$ we have 
$c^{p_{\xi+1}}_\alpha=c^{p_\xi}_\alpha$. 
\end{enumerate}
Then the sequence $\bar{p}$ has an upper bound in $\clfo$.
\item Suppose that $p\in\clfo$ and $\name{h}$ is a $\clfo$--name such that
$p\forces$``$\name{h}:\lambda\longrightarrow\bV$''. Then there is a
condition $q\in\clfo$ stronger than $p$ and such that    
\begin{enumerate}
\item[$(\otimes)$] if $\delta<\delta'$ are two successive points of $C^q$,
$u\in\pos^+_S(w^q,\gc^q,\delta)$, then the condition $q\rest_{\delta'} 
c^q_\delta(u)$ decides the value of $\name{h}\rest(\delta+1)$. 
\end{enumerate}
\end{enumerate}
\end{proposition}

\begin{proof}
(1)\quad First let us note that if $\delta\in\mathop{\triangle}\limits_{
\xi<\lambda} C^{p_\xi}$ is a limit ordinal, then $\delta\in\bigcap\limits_{
\xi<\lambda} C^{p_\xi}$ and $c^{p_{\delta+1}}_\delta=c^{p_\xi}_\delta$ for
all $\xi\geq\delta+2$ (by assumptions (b) and (c)). Now, we put $w^q=
w^{p_0}$ and $C^q=\{\delta\in\mathop{\triangle}\limits_{\xi<\lambda}
C^{p_\xi}:\delta\mbox{ is limit }\}$, and for $\delta\in C^q$ we let 
$c^q_\delta=c^{p_{\delta+1}}_\delta$ (thus defining $\gc^q=\langle 
c^q_\delta:\delta\in C^q\rangle$). It should be clear that $q=(w^q,C^q,
\gc^q)\in\clfo$ is an upper bound to $\bar{p}$.  
\medskip

\noindent (2)\quad Follows from (1) above and \ref{clbas}(4).
\end{proof}
 
\begin{definition}
\label{clnames}
We let $\name{W}$ and $\name{\eta},\name{\nu}$ be $\clfo$--names such that  
\[\forces_{\clfo}\name{W}=\bigcup\big\{w^p:p\in\Gamma_{\clfo}\big\}\]
and 
\[\begin{array}{ll}
\forces_{\clfo}&\mbox{`` }\name{\eta},\name{\nu}\in\bairel\mbox{ and if
}\langle\delta_\xi:\xi<\lambda\rangle\mbox{ is the increasing enumeration of
}\cl(\name{W}),\\
&\ \mbox{ and }\delta_\xi\leq\alpha<\delta_{\xi+1},\ \xi<\lambda,\mbox{ then
} \name{\eta}(\alpha)=\xi\mbox{ and }\name{\nu}(\alpha)=\delta_{\xi+4}
\mbox{ ''.}\end{array}\]
\end{definition}

\begin{proposition}
\label{easyW}
\begin{enumerate}
\item $\forces_\clfo$`` $\name{W}$ is an unbounded $S$--closed subset 
of $\lambda$ ''. Consequently $\forces_\clfo$`` $\name{W}\in
 \cU^{\clfo}$ ''. 
\item $\forces_\clfo\mbox{`` }\name{W},\lambda\setminus\name{W}\in
\big(\cV^\clfo\big)^+\mbox{ ''}$. 
\item $\forces_{\clfo}\big(\forall f\in\bairel{\cap}\bV\big)
  \big(\forall A\in\cV^{\clfo}\big)\big(\exists\alpha\in A\big)\big(
  f(\alpha)<\name{\nu}(\alpha)\big)$.
\end{enumerate}
\end{proposition}

\begin{proof}
(2)\quad Suppose that $p\in\clfo$ and $\name{A}_i$ (for $i<\lambda$) are
$\clfo$--names for elements of $\cV\cap\bV$. Build inductively sequences
$\langle p_i:i\leq\lambda\rangle\subseteq\clfo$ and $\langle A_i:i\leq
  \lambda\rangle\subseteq\cV$ such that 
\begin{enumerate}
\item[(a)] $\big(\forall i<j<\lambda\big)\big(p\leq p_i\leq p_j\big)$, 
\item[(b)] $p_{i+1}\forces_{\clfo}\name{A}_i=A_i$ and $i\leq\sup(w^{p_i})$
  for all $i<\lambda$,  
\item[(c)] if $\gamma<\lambda$ is limit, then  $p_\gamma$ is the natural limit of
  $\langle p_i:i<\gamma\rangle$. 
\end{enumerate}
Pick $\delta\in\mathop{\triangle}\limits_{i<\lambda} A_i\setminus S$ such
that $\delta=\sup\big(\bigcup\limits_{i<\delta} w^{p_i}\big)\in
C^{p_\delta}$ (possible by the normality of $\cV$; remember (b,c)
above). Then $p_\delta\forces\delta\in\mathop{\triangle}\limits_{i<\lambda}
\name{A}_i$. Put $\beta=\min\big(C^{p_\delta}\setminus (\delta+1)\big)$.

Let $w=c^{p_\delta}_\delta(w^{p_{\delta}})$ and $p^*=p\rest_\beta w$. Then 
$p^*\geq p$ and $p\forces \delta\in\name{W}$. 

On the other hand, since $\delta=\sup(w^{p_\delta})\notin S$, we have
$w^{p_\delta}\in\pos_S^+(w^{p_\delta},\gc^{p_\delta},\beta)$ so we may let 
$p^{**}=p\rest_\beta w^{p_\delta}$. Then $p^{**}\geq p$ and $p\forces
\delta\notin\name{W}$.  
\medskip

\noindent (3)\quad Suppose that $p\in\clfo$, $f\in\bairel$ and $\langle
\name{A}_\alpha:\alpha<\lambda\rangle$ is a sequence of $\clfo$--names for
members of $\cV\cap\bV$. By induction on $\alpha<\lambda$ construct a
sequence $\langle p_\alpha,A_\alpha:\alpha<\lambda\rangle$ such that for
each $\alpha$:
\begin{enumerate}
\item[(i)]   $p_\alpha\in\clfo$, $A_\alpha\subseteq\lambda$,
$A_\alpha\in\cV$, $p_0=p$, $p_\alpha\leq_{\clfo} p_{\alpha+1}$, and 
\item[(ii)]  if $\alpha$ is a limit ordinal, then $p_\alpha$ is the natural
limit of $\langle p_\beta:\beta<\alpha\rangle$, and  
\item[(iii)] $p_{\alpha+1}\forces_{\clfo} \name{A}_\alpha\cap (\lambda
\setminus S)=A_\alpha$. 
\end{enumerate}
Next pick a limit ordinal $\delta\in\mathop{\triangle}\limits_{\alpha<
\lambda} A_\alpha\cap (\lambda\setminus S)$ such that $(\forall\alpha<
\delta)(w^{p_\alpha}\subseteq\delta)$. Then $p_\delta\forces\delta\in
\mathop{\triangle}\limits_{\alpha<\lambda}\name{A}_\alpha$ and
$w^{p_\delta}\subseteq\delta$ is $S$--closed, so we may let
$w^q=w^{p_\delta}$, $C^q=C^{p_\delta}\setminus\big(f(\delta)+1\big)$ and
$\gc^q=\gc^{p_\delta}\rest C^q$ to get a condition $q\in\clfo$ stronger than
$p$ and such that 
\[q\forces_{\clfo}\mbox{`` }\delta\in\mathop{\triangle}\limits_{\alpha<
\lambda}\name{A}_\alpha\mbox{ and }f(\delta)<\name{\nu}(\delta)\mbox{ ''.}\] 
\end{proof}

\begin{proposition}
\label{clfoOK}
The forcing notion $\clfo$ is reasonably B--bounding over $\cU$.
\end{proposition}

\begin{proof}
By \ref{clbas}(1), $\clfo$ is $({<}\lambda)$--complete, so we have to verify 
\ref{p.1A}(5b) only. Let $p\in\clfo$ and let $\bar{\mu}=\langle
\mu_\alpha':\alpha<\lambda\rangle$, $\mu_\alpha'=\lambda$ for each
$\alpha<\lambda$. We are going to describe a strategy $\st$ for Generic in
$\tegamel(p,\clfo)$.  

In the course of a play the strategy $\st$ instructs Generic to build aside
an increasing sequence of conditions $\bar{p}^*=\langle p^*_\alpha:\alpha<
\lambda\rangle\subseteq\clfo$ such that  
\begin{enumerate}
\item[(a)] $p_0^*=p$ and $w^{p^*_\alpha}=w^p$ for all $\alpha<\lambda$, and  
\item[(b)] if $\gamma<\lambda$ is limit, then $p^*_\gamma$ is the natural
  limit of $\bar{p}^*\rest\gamma$, and 
\item[(c)] for each $\alpha<\lambda$, if $\delta\in C^{p^*_\alpha}$,
$\otp(C^{p^*_\alpha}\cap\delta)=\alpha$, then $C^{p^*_{\alpha+1}}\cap
(\delta+1)=C^{p^*_\alpha}\cap (\delta+1)$ and for every $\xi\in
C^{p^*_{\alpha+1}}\cap\delta$ we have $c^{p^*_{\alpha+1}}_\xi=
c^{p^*_\alpha}_\xi$, and 
\item[(d)] after stage $\alpha<\lambda$ of the play of $\tegamel(p,\clfo)$,
the condition $p_{\alpha+1}$ is determined (conditions $p_\alpha$ for
non-successor $\alpha<\lambda$ are determined by (a),(b) above). 
\end{enumerate}
So suppose that the players arrived to a stage $\alpha<\lambda$ of
$\tegamel(p,\clfo)$, and Generic (playing according to $\st$ so far) has
constructed aside an increasing sequence $\langle p^*_\xi:\xi\leq p^*_\alpha
\rangle$ of conditions (satisfying (a)--(d)). Let $\delta\in C^{p^*_\alpha}$
be such that $\otp(C^{p^*_\alpha}\cap\delta)=\alpha$ and let $\gamma=\min(
C^{p^*_\alpha}\setminus(\delta+1))$. Now Generic makes her move in
$\tegamel(p,\clfo)$:  
\begin{itemize}
\item $I_\alpha=\pos^+_S(w^{p^*_\alpha},\gc^{p^*_\alpha},\delta)$, and 
\item $p^\alpha_u=p^*_\alpha\rest_\gamma c^{p^*_\alpha}_\delta(u)$ for $u\in
  I_\alpha$. 
\end{itemize}
Let $\langle q^\alpha_u:u\in I_\alpha\rangle\subseteq\clfo$ be the answer of
Antigeneric, so $p^*_\alpha\rest_\gamma c^{p^*_\alpha}_\delta(u)\leq
q^\alpha_u$ for each $u\in\pos^+(w^{p^*_\alpha},\gc^{p^*_\alpha},\delta)$.
Now Generic uses \ref{clbas}(3) (with $\delta,\gamma,p^*_\alpha,q^\alpha_u$
here standing for $\alpha,\beta,p,q_u$ there) to pick a condition
$p^*_{\alpha+1}$ such that, letting $\alpha'=\min(C^{p^*_{\alpha+1}}
\setminus\gamma)$, we have  
\begin{enumerate}
\item[(e)] $p^*_\alpha\leq p^*_{\alpha+1}$, $w^{p^*_{\alpha+1}}=w^p$, $C^{
p^*_{\alpha+1}}\cap\gamma=C^{p^*_\alpha}\cap\gamma$ and $c^{p^*_{\alpha+
1}}_\xi=c^{p^*_\alpha}_\xi$ for $\xi\in C^{p^*_{\alpha+1}}\cap\delta$, and 
\item[(f)] $\bigcup\big\{w^{q^\alpha_u}:u\in I_\alpha\big\}\subseteq
  \alpha'$, and  
\item[(g)] $q^\alpha_u\leq p^*_{\alpha+1}\rest_{\alpha'}
  c^{p^*_{\alpha+1}}_\delta(u)$ for every  $u\in I_\alpha$. 
\end{enumerate}
We claim that $\st$ is a winning strategy for Generic in
$\tegamel(p,\clfo)$. So suppose that 
\[\Big\langle I_\alpha,\langle p^\alpha_u,q^\alpha_u:u\in I_\alpha\rangle:
\alpha<\lambda\Big\rangle\]  
is a play of $\tegamel(p,\clfo)$ in which Generic uses $\st$, and let 
$\bar{p}^*=\langle p^*_\alpha:\alpha<\lambda\rangle\subseteq\clfo$ be the
sequence constructed aside by Generic, so it satisfies (a)--(c) above, and
thus also the assumptions of \ref{cllimit}(1). Let $p^*$ be an upper bound to
$\bar{p}$ (which exists by \ref{cllimit}(1)). Now note that 
\[p^*\forces_\clfo\mbox{`` if }\alpha\in C^{p^*}\cap\name{W}\mbox{ and }u=
\name{W}\cap\alpha,\mbox{ then }q^\alpha_u\in\Gamma_{\clfo}\mbox{ ''}\]
and therefore
\[p^*\forces_\clfo\mbox{`` }\big(\forall\alpha\in C^{p^*}\cap \name{W}\big)
\big(\exists u\in I_\alpha\big)\big(q^\alpha_u\in\Gamma_{\clfo}\big)
\mbox{ ''.}\]
Since $p^*\forces C^{p^*}\cap\name{W}\in\cU^\clfo$ (by \ref{easyW}) we may
conclude that the condition $p^*$ witnesses that Generic won the play. 
\end{proof}

\begin{definition}
Let $\cF$ be a filter on $\lambda$ including all co=bounded subsets of
$\lambda$, $\emptyset\notin\cF$. 
\begin{enumerate}
\item We say that a family $F\subseteq\bairel$ is {\em $\cF$--dominating\/}
  whenever 
\[\big(\forall g\in \bairel\big)\big(\exists f\in F\big)
\big(\{\alpha<\lambda: g(\alpha)<f(\alpha)\}\in\cF\big).\]
\item The $\cF$--dominating number $\gd_\cF$ is the minimal size of an
  $\cF$--dominating family in $\bairel$.
\item If $\cF$ is the filter of co-bounded subsets of $\lambda$, then the
  corresponding dominating number is also denoted by $\gd_\lambda$. If $\cF$ 
  is the filter generated by club subsets of $\lambda$, then the
  corresponding dominating number is called $\gd_{\rm cl}$.   
\end{enumerate}
\end{definition}

It was shown in Cummings and Shelah \cite{CuSh:541} that
$\gd_\lambda=\gd_{\rm cl}$ (whenever $\lambda>\beth_\omega$ is regular). The
following conclusion is an interesting addition to that result. 

\begin{conclusion}
\label{conc}
It is consistent that $\lambda$ is an inaccessible cardinal and there are
two normal filters $\cU',\cU''$ on $\lambda$ such that
$\gd_{\cU'}\neq\gd_{\cU''}$.  
\end{conclusion}

\begin{proof}
Start with the universe where $\lambda,\cU,\cV,S$ are as in \ref{incon}+
\ref{extraassum} and $2^\lambda=\lambda^+$. Let $\bar{\bbQ}=\langle\bbP_\xi,
\name{\bbQ}_\xi:\xi<\lambda^{++}\rangle$ be a $\lambda$--support iteration
such that for every $\xi<\lambda^{++}$,  $\forces_{\bbP_\xi}\mbox{``
}\name{\bbQ}_\xi=\clfo$ ''. 

It follows from \ref{verB} that $\bbP_{\lambda^{++}}$ is reasonably {\bf
b}--bounding over $\cU$, and hence also $\lambda$--proper. Therefore using
\ref{clbas}(1) and \cite[Theorem A.1.10]{RoSh:777} (see also Eisworth
\cite[\S 3]{Ei03}) one can easily argue that the limit $\bbP_{\lambda^{++}}$
of the iteration satisfies the $\lambda^{++}$--cc, has a dense subset of
size $\lambda^{++}$, is strategically $({<}\lambda)$--complete and
$\lambda$--proper. Consequently, the forcing with $\bbP_{\lambda^{++}}$ does
not collapse cardinal. Also it follows from \ref{bound} that    
\[\forces_{\bbP_{\lambda^{++}}}\mbox{`` $\bairel\cap\bV$ is
  $\big(\cU\big)^{\bbP_{\lambda^{++}}}$--dominating in $\bairel$ ''}\]
and it follows from \ref{easyW}(3) that for each $\xi<\lambda^{++}$
\[\forces_{\bbP_{\lambda^{++}}}\mbox{`` $\bairel\cap\bV^{\bbP_\xi}$ is not  
  $\big(\cV\big)^{\bbP_{\lambda^{++}}}$--dominating in $\bairel$ ''}\]
Therefore we may easily conclude that 
\[\begin{array}{ll}
\forces_{\bbP_{\lambda^{++}}}&\mbox{`` if }\cU'=\big(\cU
\big)^{\bbP_{\lambda^{++}}},\ \cU''=\big(\cV\big)^{\bbP_{\lambda^{++}}}
\mbox{ then}\\
&\quad\gb_\lambda=\gd_{\cU'}=\lambda^+<2^\lambda=\lambda^{++}=\gd_{\cU''}=
\gd_{\rm cl}=\gd_\lambda\mbox{ ''.}
  \end{array}\]
\end{proof}

\[*\qquad *\qquad *\]

\begin{definition}
\label{3.2}
We define a forcing notion $\tefo$ as follows.\\
{\bf A condition in $\tefo$} is a triple $p=(w^p,C^p,\gc^p)$ such that 
\begin{enumerate}
\item[(i)]    $C^p\in \cU$, $w^p\subseteq\min(C^p)$,  
\item[(ii)]   $\gc^p=\langle c^p_\alpha:\alpha\in C^p\rangle$ is a
  $C^p$--extending sequence.
\end{enumerate}
{\bf The order $\leq_{\tefo}=\leq$ of $\tefo$} is given by \\
$p\leq_{\tefo} q$\qquad if and only if \\
\begin{enumerate}
\item[(a)] $C^q\subseteq C^p$ and $w^q\in\pos^+(w^p,\gc^p,\min(C^q))$ and  
\item[(b)] if $\alpha_0,\alpha_1\in C^q$, $\alpha_0<\alpha_1=\min(C^q
\setminus(\alpha_0+1))$ and $u\in\pos^+(w^q,\gc^q,\alpha_0)$, then
$c^q_{\alpha_0}(u)\in \pos(u,\gc^p,\alpha_1)$.
\end{enumerate}
For $p\in\tefo$, $\alpha\in C^p$ and $u\in\pos^+(w^p,\gc^p,\alpha)$ we let  
$p\rest_\alpha u\stackrel{\rm def}{=}(u,C^p\setminus\alpha,\gc^p\rest(C^p
\setminus\alpha))$.
\end{definition}

\begin{proposition}
\label{3.4}
\begin{enumerate}
\item $\tefo$ is a $\lambda$--complete forcing notion of cardinality
  $2^\lambda$. 
\item If $p\in\tefo$ and $\alpha\in C^p$, then
\begin{itemize}
\item for each $u\in\pos^+(w^p,\gc^p,\alpha)$, $p\rest_\alpha u\in\tefo$ is
  a condition stronger than $p$, and    
\item the family $\{p\rest_\alpha u: u\in\pos^+(w^p,\gc^p,\alpha)\}$ is
  pre-dense above $p$. 
\end{itemize}
\item Let $p\in\tefo$ and $\alpha<\beta$ be two successive members of
$C^p$. Suppose that for each $u\in\pos^+(w^p,\gc^p,\alpha)$ we are given a
condition $q_u\in\tefo$ such that $p\rest_\beta c^p_\alpha(u)\leq q_u$. 
Then there is a condition $q\in\tefo$ such that letting $\alpha'=\min(
C^q\setminus\beta)$ we have
\begin{enumerate}
\item[(a)] $p\leq q$, $w^q=w^p$, $C^q\cap\beta=C^p\cap\beta$ and
  $c^q_\delta=c^p_\delta$ for $\delta\in C^q\cap\alpha$, and 
\item[(b)] $\bigcup\big\{w^{q_u}:u\in\pos^+(w^p,\gc^p,\alpha)\big\}\subseteq 
  \alpha'$, and 
\item[(c)] $q_u\leq q\rest_{\alpha'} c^q_\alpha(u)$ for every
  $u\in\pos^+(w^p,\gc^p,\alpha)$. 
\end{enumerate}
\item Assume that $p\in\tefo$, $\alpha\in C^p$ and $\name{\tau}$ is a
$\tefo$--name such that $p\forces$``$\name{\tau}\in\bV$''. Then there is a
  condition $q\in\tefo$ stronger than $p$ and such that   
\begin{enumerate}
\item[(a)] $w^q=w^p$, $\alpha\in C^q$ and $C^q\cap\alpha=C^p\cap\alpha$,
  and 
\item[(b)] if $u\in\pos^+(w^q,\gc^q,\alpha)$ and $\gamma=\min(C^q\setminus 
(\alpha+1))$, then the condition $q\rest_\gamma c^q(u)$ forces a value to
$\name{\tau}$. 
\end{enumerate}
\end{enumerate}
\end{proposition}

\begin{proof}  
Fully parallel to \ref{clbas}.
\end{proof}

\begin{definition}
{\em The natural limit\/} of an $\leq_{\tefo}$--increasing sequence
$\bar{p}=\langle p_\xi:\xi<\lambda\rangle\subseteq\tefo$ (where
$\gamma<\lambda$ is a limit ordinal) is the condition $q=(w^q,C^q,\gc^q)$
defined as follows: 
\begin{itemize}
\item $w^q=\bigcup\limits_{\xi<\gamma}w^{p_\xi}$, $C^q=\bigcap\limits_{\xi<
\gamma} C^{p_\xi}$ and 
\item $\gc^q=\langle c^q_\delta:\delta\in C^q\rangle$ is such that for
  $\delta\in C^q$ and $u\subseteq\delta$ we have $c^q_\delta(u)=
  \bigcup\limits_{\xi<\gamma}c^{p_\xi}_\delta(u)$. 
\end{itemize}
\end{definition}

\begin{proposition}
\label{limit}
\begin{enumerate}
\item Suppose $\bar{p}=\langle p_\xi:\xi<\lambda\rangle$ is a
$\leq_{\tefo}$--increasing sequence of conditions from $\tefo$ such that 
\begin{enumerate}
\item[(a)] $w^{p_\xi}=w^{p_0}$ for all $\xi<\lambda$, and 
\item[(b)] if $\gamma<\lambda$ is limit, then $p_\gamma$ is the natural
  limit of $\bar{p}\rest\gamma$, and 
\item[(c)] for each $\xi<\lambda$, if $\delta\in C^{p_\xi}$, $\otp(C^{p_\xi}
\cap\delta)=\xi$, then $C^{p_{\xi+1}}\cap (\delta+1)=C^{p_\xi}\cap
(\delta+1)$ and for every $\alpha\in C^{p_{\xi+1}}\cap\delta$ we have 
$c^{p_{\xi+1}}_\alpha=c^{p_\xi}_\alpha$. 
\end{enumerate}
Then the sequence $\bar{p}$ has an upper bound in $\tefo$.
\item Suppose that $p\in\tefo$ and $\name{h}$ is a $\tefo$--name such that
$p\forces$``$\name{h}:\lambda\longrightarrow\bV$''. Then there is a
condition $q\in\tefo$ stronger than $p$ and such that    
\begin{enumerate}
\item[$(\otimes)$] if $\delta<\delta'$ are two successive points of $C^q$,
$u\in\pos(w^q,\gc^q,\delta)$, then the condition $q\rest_{\delta'}
c^q_\delta(u)$ decides the value of $\name{h}\rest(\delta+1)$. 
\end{enumerate}
\end{enumerate}
\end{proposition}

\begin{proof}
Fully parallel to \ref{cllimit}.
\end{proof}
 
\begin{definition}
\label{names}
We let $\name{W}$ and $\name{\eta},\name{\nu}$ be $\tefo$--names such that  
\[\forces_{\tefo}\name{W}=\bigcup\big\{w^p:p\in\Gamma_{\tefo}\big\}\]
and 
\[\begin{array}{ll}
\forces_{\tefo}&\mbox{`` }\name{\eta},\name{\nu}\in\bairel\mbox{ and if
}\langle\delta_\xi:\xi<\lambda\rangle\mbox{ is the increasing enumeration of
}\cl(\name{W}),\\
&\ \mbox{ and }\delta_\xi\leq\alpha<\delta_{\xi+1},\ \xi<\lambda,\mbox{ then
} \name{\eta}(\alpha)=\xi\mbox{ and }\name{\nu}(\alpha)=\delta_{\xi+4}
\mbox{ ''.}\end{array}\]
\end{definition}

Note that if $p\in\tefo$, then 
\[p\forces_{\tefo}\mbox{`` }\name{W}\subseteq\bigcup\big\{[\alpha_0,
\alpha_1):\alpha_0,\alpha_1\in C^p\ \&\ \alpha_1=\min\big(C^p\setminus(
\alpha_0+1)\big)\big\}\mbox{ ''}\]
and 
\[\begin{array}{ll}
p\forces_{\tefo}&\mbox{`` }\big\{\alpha\in C^p:\big[\alpha,\min(C^p
\setminus(\alpha+1))\big)\cap\name{W}\neq\emptyset\big\}, \\
&\ \ \big\{\alpha\in C^p:\big[\alpha,\min(C^p\setminus
(\alpha+1))\big)\cap\name{W}=\emptyset\big\}\in\big(\cU^\tefo\big)^+
\mbox{ ''}.\end{array}\] 

\begin{proposition}
\label{3.1} 
$\forces_{\tefo}\big(\forall f\!\in\!\bairel{\cap}\bV\big)
\big(\forall A\!\in\!\cU^{\tefo}\big)\big(\exists\alpha\!\in\! A\big)
\big(f(\alpha)<\name{\nu}(\alpha)\big)$.
\end{proposition}

\begin{proof}
Fully parallel to \ref{easyW}.
\end{proof}

\begin{proposition}
\label{tefoOK}
The forcing notion $\tefo$ is reasonably C--bounding over $\cU$.
\end{proposition}

\begin{proof}
Fully parallel to \ref{clfoOK}.
\end{proof}

The following problem is a particular case of \ref{prob}(1).

\begin{problem}
Are $\lambda$--support iterations of $\tefo$ $\lambda$--proper? 
\end{problem}

\end{document}